\documentclass[12pt,a4paper]{amsart}
\usepackage{a4wide}
\usepackage{amsmath, amssymb, amsfonts, enumerate}

 \theoremstyle{definition}

\numberwithin{equation}{section}

\newcommand {\N}{\mathbf{N}} 
\newcommand {\Z}{\mathbf{Z}} 
\newcommand {\R}{\mathbf{R}} 
\newcommand {\C}{\mathbf{C}} 

\begin{document}
\title[Amenability and paradoxical decompositions]{Amenability and paradoxical decompositions for pseudogroups and for discrete metric spaces}
\author[T. Ceccherini-Silberstein]{Tullio Ceccherini-Silberstein}
\address{Dipartimento di Ingegneria, Universit\`a del Sannio, C.so
Garibaldi 107, 82100 Benevento, Italy}
\email{tceccher@mat.uniroma3.it}
\author[R.I. Grigorchuk]{Rostislav I. Grigorchuk}
\address{Mathematics Department,
Texas A \& M University, College Station, TX 77843-3368, USA}
\email{grigorch@math.tamu.edu}
\author[P. de la Harpe]{Pierre de la Harpe}
\address{Section de Math\'ematiques,  Section de math\'ematiques, 2--4 rue du Li\`evre C.P. 64. CH 1211 Gen\`eve 4. Switzerland}
\email{Pierre.delaHarpe@unige.ch}
\thanks{The  authors acknowledge support from the Fonds National Suisse de la Recherche Scientifique}

\date{\today}
\begin{abstract}
This is an expostion of various aspects of amenability and paradoxical decompositions for groups, group actions and metric spaces. First, we review the formalism of pseudogroups, which is well adapted to stating the alternative of  Tarski,  according to which a pseudogroup without invariant mean gives rise to paradoxical decompositions, and to defining a F\o lner condition. Using a  Hall-Rado Theorem on matchings in graphs, we show then for pseudogroups that existence of an invariant mean is equivalent to the F\o lner condition; in the case of the pseudogroup of bounded perturbations of the identity on a locally finite metric space,   these conditions are moreover equivalent to the negation of the  Gromov's so-called doubling condition, to isoperimetric conditions, to Kesten's spectral condition for  related simple random walks, and to various other conditions. We define also the  minimal Tarski number of  paradoxical decompositions associated to a non-amenable group action (an integer $\ge 4$),  and we indicate numerical estimates (Sections II.4 and IV.2). The final chapter explores for metric spaces the notion of supramenability, due for groups to Rosenblatt.
\end{abstract}
\maketitle

\section*{{\bf I. Introduction}}
The present exposition shows various aspects of amenability and  non-amenability. Our initial motivation comes from a note on the {\it Banach-Tarski paradox} where Deuber, Simonovitz and S\'os indicate one kind of paradoxical decomposition for metric spaces, in relation with what they call  an ``exponential growth'' property \cite{DeSS}. Our first purpose is  to revisit their work which, in our view, relates paradoxical decompositions  with {\it amenability} rather than with growth 
(see in particular Observation  33 below). 
\par     
For this, we recall in Chapter II the formalism of set-theoretical pseudogroups which is well adapted to showing the many aspects of amenability: existence of {\it invariant finitely additive measures}, absence of {\it paradoxical decompositions}, existence of {\it F\o lner sets} and {\it isoperimetric estimates}. We also state one version of the  basic {\it Tarski alternative}: a pseudogroup is either amenable or paradoxical. 
\par     
In Chapter III, we specialize the discussion to  metric spaces and pseudogroups of bounded perturbations of the identity; metric spaces, there, are locally finite 
(except at the very end of the chapter). On one hand, this is an interesting class, with many examples given by finitely generated groups. On the other hand, it provides a convenient setting for proving F\o lner characterization as stated in Chapter II. We discuss also the {\it Kesten characterization in terms of simple random walks}. 
\par     
For a group $G$ which is not amenable, we estimate in Chapter IV the {\it Tarski number} ${\mathcal T}(G) \in \{4, 5, \ldots , \infty \}$, 
which indicates the minimal number of pieces involved in a paradoxical decomposition of $G$. It is known that ${\mathcal T}(G) = 4$ if and only if $G$ has a subgroup which is free non abelian. We show that one has $5 \le {\mathcal T}(G) \le 34$   [respectively $6 \le {\mathcal T}(G) \le 34$] for some torsion-free groups [resp. for some torsion groups]  constructed by Ol'shanskii \cite{Ol1}, and $6 \le {\mathcal T}\big(B(m,n)\big) \le 14$ for $B(m,n)$ a {\it Burnside group} on $m \ge 2$ generators of odd exponent $n \ge 665$ \cite{Ady2}.  
\par     
Building upon the seminal 1929 paper by von Neumann \cite{NeuJ}, Rosenblatt has defined for groups a notion of {\it supramenability.} He has shown that supramenable groups include those of subexponential growth, and it is not  known whether there are others. In Chapter V, we investigate supramenability for pseudogroups and for locally finite metric spaces; in particular, we describe a simple example of a graph which is both supramenable and of superexponential growth.   
\par   
We are grateful to Joseph Dodziuk, Vadim Kaimanovich, Alain Valette and Wolfgang Woess for useful discussions and bibliographical informations, as well as to Laurent Bartholdi for Presentation 12, Example 74, and his critical reading of a preliminary version of this work.\footnote{ 
This post on arXiv is the published version (Proc. Steklov Inst. Math.  {\bf 224} (1999), 57--97) with the following changes: (i) the caution following Definition 28, (ii) the addition of a missing hypothesis in Proposition 38 
(we are grateful to Volker Diekert for having pointed out this omission to us), (iii) the updating of some references, and
(iv) the correction of a few minor typos.

Moreover, we have collected comments on several items in a new Chapter VI, after the first list of references.}

\section*{{\bf II. Amenable pseudogroups}} 
\section*{II.1. Pseudogroups} 

\noindent
{\bf 1. Definition.} In the present set-theoretical context, a {\it pseudogroup} ${\mathcal G}$ of transformations of a set $X$ is a set of bijections $\gamma \colon S \to T$ between subsets $S,T$ of $X$ which satisfies the following conditions (as listed, e.g., in \cite{HS1}):
\begin{enumerate}[{\rm (i)}]

\item the identity $X \to X$ is in ${\mathcal G}$,
 
\item if $\gamma \colon S \to T$ is in ${\mathcal G}$, so is the inverse $\gamma^{-1}\colon T \to S$, 

\item if $\gamma \colon S \to T$ and $\delta \colon T \to U$ are in ${\mathcal G}$,   so is their composition $\delta \gamma \colon S \to U$, 

\item if  $\gamma \colon S \to T$ is in ${\mathcal G}$ and if $S'$ is a subset of $S$, the restriction $\gamma \vert S' \colon S' \to \gamma(S')$ is in  ${\mathcal G}$, 

\item if $\gamma \colon S \to T$ is a bijection between two subsets $S,T$ of $X$  and 
if there exists a {\it finite} partition $S = \sqcup_{1 \le j \le n} S_j$ with $\gamma \vert S_j$ in ${\mathcal G}$ for $j \in \{1, \hdots, n\}$,  then $\gamma$ is in ${\mathcal G}$   (where $\sqcup$ denotes a disjoint union).
\end{enumerate}
Property (v) expresses the fact that ${\mathcal G}$ is closed with respect to {\it finite gluing up}; together with (iv), they express the fact that, for a bijection $\gamma$, being in ${\mathcal G}$ is in some sense a {\it local} condition. 
\par 

For $\gamma \colon S \to T$ in ${\mathcal G}$, we write also $\alpha (\gamma)$
for the domain $S$ of $\gamma$ and $\omega (\gamma)$ for its range $T$. For
``a pseudogroup ${\mathcal G}$ of transformations of a set $X$'' , we
write shortly
``a pseudogroup $({\mathcal G} , X)$'', or even 
``a pseudogroup ${\mathcal G}$''.

\bigskip

\noindent
{\bf 2. Examples.} (i) Any {\it action of a group} $G$ on a set
$X$ generates a pseudogroup ${\mathcal G}_{G,X}$. More precisely, a bijection 
$\gamma \colon S \to T$ is in ${\mathcal G}_{G,X}$ if there exists a finite
partition $S = \sqcup_{1 \le j \le n} S_j$
and elements $g_1, \hdots, g_n \in G$ such that
$\gamma(x) = g_j(x)$ for all $x \in S_j, j \in \{1, \hdots, n\}$.
If there exists such a $\gamma$, the subsets $S,T$ of $X$ are sometimes said 
to be $G$-equidecomposable (or ``endlich zerlegungsgleich''\, in
\cite{NeuJ}). \par
   In case $G = X$ acts on itself by left multiplications, we write
${\mathcal G} _G$ instead of ${\mathcal G}_{G,G}$.\smallskip

   (ii) {\it Piecewise isometries} of a metric space $X$ constitute a
pseudogroup ${\mathcal P}i{\mathcal I}s(X)$, generated (in the obvious way) by
the partial isometries between subsets of $X$. Observe that it may be much
larger than the pseudogroup associated as in the previous example with the
group of isometries of $X$;  see for example the metric space obtained from
the real line by gluing two hairs of different length at two distinct points
of the line.\smallskip

   (iii) For a metric space $X$, the pseudogroup  ${\mathcal W}(X)$ of  
{\it bounded perturbations of the identity} consists of bijections  
$\gamma \colon S \to T$ such that $\sup_{x \in S}d(\gamma (x) , x) < \infty$. 
In agreement with the main example of \cite{DeSS}, we like to call ${\mathcal W}(X)$ the pseudogroup of
{\it wobbling bijections}; the notion seems to come from the important work
by Laczkovich \cite{Lacz}. See also Item $0.5.C_1''$ in \cite{Gro3}.\smallskip

   (iv) Given a pseudogroup ${\mathcal G}$ of transformations of a set $X$ and a  
{\it subset} $A$ of $X$, the set of bijections $\gamma \in {\mathcal G}$ with
$\alpha (\gamma) \subset A$ and $\omega(\gamma) \subset A$ constitute a
pseudogroup of transformations of $A$, denoted below by ${\mathcal G}_{(A)}$.\smallskip

   (v) From a pseudogroup $({\mathcal G} , X)$ and an integer $k \ge 1$, one 
obtains a pseudogroup ${\mathcal G}_k$ of transformations of the direct product
$X_k$ of $X$ and $\{1, \hdots, k\}$, generated by the bijections of the 
form 
$$
\left\{
\aligned
   S \times \{j\} \quad &\longrightarrow \quad T \times \{j'\} \\
    (x,j) \quad &\longmapsto \quad (\gamma(x) , j')
\endaligned
\right.
$$
where $\gamma \colon S \to T$ is in ${\mathcal G}$ and $1 \le j , j' \le k$.
\bigskip

\noindent
{\bf 3. Remarks.} The above notion of pseudogroup of
transformations is strongly motivated by the study of Banach-Tarski
paradoxes, as shown by the first three observations below.
\smallskip

   (i) The very definition of a paradoxical decomposition with respect to a
group action involves the associated pseudogroup as in Example 2.(i).\smallskip

   (ii) Pseudogroups are easily restricted on subsets as in Example 2.(iv).
This is important for the study of supramenability (see Chapter V below).
\smallskip

   (iii) Pseudogroups are easily induced on oversets, as in Example
2.(v). This is useful in the setting of a  pseudogroup constituted by 
bijections  with domains and range required to be in a given
algebra  (or $\sigma$-algebra) of subsets of $X$ (for example the measurable
sets of a measure space), and in corresponding variations
on the Tarski alternative \cite{HS1}. \smallskip

   (iv) For a pseudogroup $({\mathcal G} ,X)$, the set
$$
   {\mathcal R} = \left\{(x,y) \in X \times X \, \big\vert \, 
    \text{there exists $\gamma \in {\mathcal G}$ such that $x \in \alpha(\gamma)$
           and $y = \gamma (x)$}\right\}
$$
is an {\it equivalence relation.} A natural problem is to study the existence of
measures $\mu$ on $X$ such that, for each measurable subset $A$ of $X$ of
measure zero, the saturated set $\{x \in A \, \vert \,  \mbox{ there exists }  a \in A \text{ with } (x,a) \in {\mathcal R}\}$ has also measure
zero, see \cite{CoFW}, \cite{Kai2}, \cite{Kai3}. \smallskip

   (v) In a topological context, Conditions (iv) and (v) in Definition 1
are usually replaced by a condition involving restrictions to {\it open
subsets}; see \cite{Sac} and page 1 of \cite{KoNo}.\smallskip

   (vi) Consider a metric space $X$, the pseudogroup ${\mathcal W}(X)$ of
Example 2.(iii), and a subspace $A$ of $X$. It is then remarkable
(though straightforward to check) that the pseudogroup ${\mathcal W}(A)$ coincides
with the restriction of ${\mathcal W}(X)$ to $A$ in the sense of Example 2.(iv).

\section*{II.2. Amenability and paradoxical decompositions - the Tarski alternative}

\noindent
Let  $({\mathcal G} ,X)$ be a pseudogroup. We denote by ${\mathcal P}(X)$ the set of
all subsets of $X$.
\medskip

\noindent
{\bf 4. Definitions.} A {\it ${\mathcal G}$-invariant mean} on $X$ is a
mapping $\mu \colon {\mathcal P}(X) \to [0,1]$ which is \medskip

$(fa)$  \, \, finitely additive: 
$\mu(S_1 \cup S_2) \, = \, \mu(S_1) + \mu(S_2)$
for $S_1,S_2 \subset X$ with $S_1 \cap S_2 = \emptyset$, \par

$(in)$  \, \, invariant: 
$\mu\big( \omega(\gamma) \big) \, = \, \mu\big( \alpha(\gamma) \big)$
for all $\gamma \in {\mathcal G}$, \par

$(no)$  \, \, normalized: $\mu(X) \, = \, 1$. \medskip

\noindent More generally, for $A \subset X$, a 
{\it ${\mathcal G}$-invariant mean on $X$ normalized on $A$} is a mapping
$\mu \colon {\mathcal P}(X) \to [0,\infty]$ which satisfies Conditions $(fa)$ and $(in)$
above, as well as \medskip

$(no')$  \, \, $\mu(A) \, = \, 1$. 
\bigskip

The pseudogroup ${\mathcal G}$ is {\it amenable} if there exists a ${\mathcal G}$-invariant 
mean on $X$, and the triple $({\mathcal G}, X , A)$ is {\it
amenable} if there exists a ${\mathcal G}$-invariant mean on $X$ normalized on
$A$. These notions are essentially due to von Neumann \cite{NeuJ}.  \medskip

\noindent   {\bf 5. Definition.}
A {\it paradoxical ${\mathcal G}$-decomposition of $X$} is a partition
$X = X_1 \sqcup X_2$ such that there exist $\gamma_j \in {\mathcal G}$
with $\alpha (\gamma_j) = X_j$ and $\omega (\gamma_j) = X$ $(j = 1,2)$. \par
A pseudogroup $({\mathcal G} , X)$ is {\it paradoxical} if it has a paradoxical
${\mathcal G}$-decomposition, or equivalently (because of Theorem 7 below) if
it is not amenable.
\medskip

\noindent {\bf 6. Remarks.} 
(i) {\it There cannot exist such paradoxical  ${\mathcal G}$-decomposition if ${\mathcal G}$ is amenable.} 

This is obvious, because (with the notation of Definitions 4 and 5) one 
cannot have $1 = \mu (X) = \mu(X_1) + \mu(X_2) = 2$ ! \par

It is remarkable that there is no further obstruction, as Theorem 7
shows.\smallskip

(ii) Let ${\mathcal G}$ and ${\mathcal H}$ be two pseudogroups of transformations 
of the same set $X,$ with ${\mathcal G} \subset {\mathcal H}$. If ${\mathcal H}$ is amenable,
then so is ${\mathcal G};$ if ${\mathcal G}$ is paradoxical, then so is ${\mathcal H}$.
This will be used for example in the proof of Theorem 25 (Item 36). \smallskip

(iii) In short-hand, Definition 5 reads $2[X] \overset{!!}{=} [X]$. It has
variations in the literature; for example, one may ask  
$(n+1)[X] \overset{!!}{\le} n[X]$, or more precisely: 
\smallskip

there exists an integer $n \ge 1$ and elements 
$\gamma_1 , \hdots , \gamma_N \in {\mathcal G}$ such that 
\par

$\left\vert \left\{ 
j \in \{1 , \hdots , N\} \mid x \in \alpha(\gamma_j) 
\right\} \right\vert 
\ge n+1$ for all $x \in X$, namely 
$\sum_{j=1}^k [ \alpha(\gamma_j)] \ge (n+1)[X]$, 
\indent
and 
\par

$\big\vert \left\{ j \in \{1 , \hdots , N\} \mid x \in \omega(\gamma_j) 
          \right\}^{}_{} \big\vert \le n$ for all $x \in X$, 
namely
$\sum_{j=1}^k [ \omega(\gamma_j) ] \le n[X]$. 
\smallskip

\noindent
Then Remark (i) still holds for the same obvious kind of reason. Indeed, the
variation is equivalent to Definition 5, as can be seen either with manipulations 
\`a la Cantor-Bernstein  (see for example \cite{HS1}) or as a consequence of the
following theorem.
\bigskip

\noindent
{\bf 7. Theorem (Tarski alternative).} {\it Let ${\mathcal G}$ be a pseudogroup of
transformations of a set $X$. Exactly one of the following holds: \smallskip

   - either ${\mathcal G}$ is amenable, \par
   - or there exists a paradoxical ${\mathcal G}$-decomposition of $X$. \medskip

\noindent
Let moreover $A$ be a non-empty subset of $X$ and let ${\mathcal G}_{(A)}$ be
the pseudogroup obtained by restriction of ${\mathcal G}$, as in Example 2.(iv).
Exactly one of the following holds: \smallskip

   - either there exists a ${\mathcal G}$-invariant mean on $X$ normalized on $A$, 
\par
   - or there exists a paradoxical ${\mathcal G}_{(A)}$-decomposition of $A$.}
\bigskip

   The theorem originates in Tarski's work: see \cite{Tar3}, as well as 
earlier papers by Tarski (\cite{Tar1}, \cite{Tar2}). \par
    One proof for pseudogroups has been written up in
\cite{HS1}. Its starting point is an application of the Hahn-Banach
theorem, to the Banach space $\ell^{\infty} (X)$ of bounded real-valued 
functions on $X$,  to the subspace 
$d^{\infty}(X)$ of finite linear combinations of functions of the form 
$\chi \big( \omega (\gamma) \big) - \chi \big( \alpha (\gamma) \big)$
for some $\gamma \in {\mathcal G}$ (where $\chi(A)$ denotes the characteristic
function of $A$), and to the open cone ${\mathcal C}$ of functions 
$F \in \ell^{\infty}(X)$ such that $\inf_{x \in X} F (x) > 0;$
one has to observe that ${\mathcal G}$ has an invariant mean if and only if      
$d^{\infty}(X) \cap {\mathcal C}= \emptyset$. This proof 
uses also ideas of Banach, Cantor-Bernstein, Hausdorff, K\"onig, Kuratowski
and von Neumann. \par

   We give here another proof, based on what we call the   
Hall-Rado theorem (Theorem 35), which is essentially 
the ``K\"onig theorem'' of \cite{Wag}. 
More precisely, the first statement of Theorem 7 is a straightforward
consequence of Theorem 25 and Theorem 32, and the second statement follows
(see the sketch below). \par

   Much more complete information on all this can be found
in Wagon's book (see \cite{Wag}, in particular Corollary 9.2 on page 128).
Important more recent work in this area include \cite{DouF}. \medskip

   Let us  sketch  the proof of the second statement  
of the theorem. Assume that the pseudogroup
${\mathcal G}_{(A)}$ is not paradoxical, so that,
by the first statement, there exists a
${\mathcal G}_{(A)}$-invariant mean $\mu_A \colon {\mathcal P}(A) \to [0,1]$.
Define then a mapping $\mu \colon {\mathcal P}(X) \to [0,\infty]$ as follows; for a
subset $Y$ of $X$, if there exists a partition 
$Y = \sqcup_{1 \le j \le n} Y_j$ and elements 
$\gamma_1 \colon Y_1 \to B_1 , \hdots , \gamma_n  \colon Y_n \to B_n$
in ${\mathcal G}$ with $B_1 , \hdots , B_n \subset A$, then set
$\mu (Y) = \sum_{j=1}^n \mu_A (B_j);$ otherwise, set
$\mu (Y) = \infty$. Then one checks that $\mu$ is well defined
and that it is a  ${\mathcal G}$-invariant mean on $X$ normalized on $A$.
\bigskip

\noindent {\bf 8. Remark.} A famous theorem of E. Hopf can be
expressed  very much like Tarski's alternative. \par

Let $T \colon X \longrightarrow X$ be an ergodic non-singular transformation of 
a finite probability space $(X, {\mathcal B}, m)$, with $m$ non-atomic.
Let $[[T]]$ denote the set of all $1$-$1$ non-singular transformations
$\phi \colon U \to V$ such that $\phi(x)$ belongs to the
$T$-orbit of $x$ for all $x \in U$ (with $U,V \in {\mathcal B}$);
this $[[T]]$ is the {\it full groupoid of $T$} of Katznelson and Weiss 
\cite[page 324]{KaWe}. 
For two measurable subsets $A,B$ of $X$,
say that $A$  {\it is dominated by} $B$, and write $A \prec B$, 
if there exists a measurable subset $B'$ of $B$ 
with $m(B \setminus B') > 0$ and a bijective transformation
$\phi \colon A \longrightarrow B'$ in $[[T]]$. \medskip

\noindent
{\bf Hopf alternative.} {\it (i) In the situation above, exactly one of the
following holds: \smallskip

   - there exists a $T$-invariant probability measure  on 
$(X, {\mathcal B})$ equivalent to $m$,  \par

   - one has $X \prec X$. \smallskip

\noindent
   (ii) Also, exactly one of the
following holds: \smallskip 

   - there exists a $T$-invariant infinite measure  on 
$(X, {\mathcal B})$ equivalent to $m$, \par

   - one has $X \prec X$, and there exists $A \in {\mathcal B}$ with $m(A) > 0$ such
that $A$ is not dominated by $A$.}
\medskip

   In other words, $(i)$ says that there is a finite invariant measure in the
measure class $m$ if and only if $X$ itself is {\it not} 
``Hopf-compressible'', and (ii) that there is an  infinite invariant
measure in the measure class $m$ if and only if $X$ is Hopf-compressible and 
{\it some} measurable subset of $X$ of positive measure is not Hopf-compressible
\cite{Weis}.
\par

   If there exists a $T$-invariant probability measure [respectively infinite
measure]  on  $(X, {\mathcal B})$ equivalent to $m$, then $T$ is said to be of 
{\it type $II_1$}  [resp. of {\it Type $II_{\infty}$}].

 \section*{II.3. The case of groups}

\noindent
For any group $G$, we consider first the pseudogroup 
${\mathcal G}_G$ which is associated with the action of $G$ on itself on the left,  
as in Example 2.(i).
\smallskip

Let now $G$ be a group generated by a finite set $S$. 
Let $\ell_S \colon G \to \N$ denote the corresponding word length function;
thus $\ell_S$ associates with $g \in G$ the smallest integer $n \ge 0$
for which there exist $s_1 , \hdots , s_n \in S \cup S^{-1}$ with 
$g = s_1 \hdots s_n$. Let $d_L$ and $d_R$ denote respectively the left
and right invariant metrics on $G$ defined by
$$
\aligned
      d_L(x,y) \, &= \, \ell_S\left( x^{-1}y \right) \\
      d_R(x,y) \, &= \, \ell_S\left( xy^{-1} \right)
\endaligned
$$
for all $x,y \in G$. \par

   Besides ${\mathcal G}_G$, we consider also the pseudogroup 
${\mathcal P}i {\mathcal I}s (G)$ of piecewise isometries of the metric space
$(G,d_L)$, as in Example 2.(ii),  as well as the pseudogroup ${\mathcal W}(G)$     
of bounded perturbations of the identity of the metric space $(G, d_R)$,
as in Example 2.(iii). It is easy to check that the pseudogroup
${\mathcal W}(G)$ does not depend on the choice of $S$.\bigskip

\noindent
{\bf 9. Observation} {\it With the notation above, one has
${\mathcal G}_G = {\mathcal W}(G)$ for any finitely generated group $G$.}

\medskip

\noindent
{\it Proof.} It is obvious that ${\mathcal G}_G \subset {\mathcal W}(G)$. 
Conversely, let $\gamma \colon U \to V$ be in ${\mathcal W}(G)$. Set
$$
\aligned
   k \, &= \, \sup_{x \in U} d_R(\gamma (x) , x) \\
   B \, &= \, \left\{ \, g \in G \, \vert \, \ell_S(g) \le k \, \right\} 
\endaligned
$$
and observe that $B$ is a finite subset of $G$.
For each $g \in B$, set 
$$
     U_g \, = \,  \{ \, x \in U \, \vert \, \gamma (x) = gx \, \}.
$$
One has $U = \sqcup _{g \in B} U_g$ and
$\gamma (x) = gx$ for all $x \in U_g$. Hence $\gamma \in {\mathcal G}_G$. \hfill \qed

\bigskip

It is clear that ${\mathcal G}_G \subset {\mathcal P}i {\mathcal I}s (G)$. 
It is also clear that ${\mathcal G}_G \ne {\mathcal P}i {\mathcal I}s (G)$ in general
(example: for $G = \Z$ generated by $\{1\}$, the isometry
$n \mapsto -n$ is not in ${\mathcal G}_{\Z}$).

\bigskip

\noindent {\bf 10. Definition.} A group $G$ is {\it amenable} if the
pseudogroup ${\mathcal G}_G$ is amenable. 
\par
   If $G$ is finitely generated, the previous observation shows that one may
equivalently define $G$ to be amenable if the pseudogroup ${\mathcal W}(G)$ is
amenable. \bigskip

\noindent {\bf 11. On the class of amenable groups.}
   Amenability may be viewed as a finiteness condition. One of the main
problems is to understand various classes  of amenable groups, for example
those which are finitely generated or finitely presented. (Recall that a group
is amenable if and only if all its finitely generated subgroups are amenable;
see Theorem 1.2.7 in \cite{Gre1} and Observation 19 below.) \par

   The following question, implicit in \cite{NeuJ}, was formulated explicitely
by Day, at the end of Section 4 in \cite{Day1}: does every non-amenable group 
contain a free group on $2$ generators? As much as we know and despite several
misleading allusions in the literature to some 
``von Neumann conjecture'', von Neumann himself  has {\it never
conjectured} that a non-amenable group should contain a non-abelian free 
subgroup! \par

   Day's question was answered negatively by A. Yu. Ol'shanskii \cite{Ol1}, 
Adyan \cite{Ady2} and Gromov \cite[Corollary 5.6.D]{Gro2}; the first two use
cogrowth criteria (see Item 52 below) and Gromov uses  Property (T). 
   For infinite groups, this {\it Property (T) of Kazhdan} \cite{Kaz} is  
(among other things) a strong form of non-amenability: see \cite{Sch} and
\cite{CoWe}. However, when restricted to the class of {\it linear groups}
(i.e. of groups which have faithful finite-dimensional linear representations),
Day's question can be answered positively: this follows from an important
result due to Tits \cite{Tit}.
\par

   M. Day has defined the class EG of ``elementary amenable groups'',
which is the smallest class of groups which contains finite groups and abelian
groups, and which is closed under the four operations of $(i)$ taking
subgroups, $(ii)$ forming factor groups, $(iii)$ group extensions and
$(iv)$ upwards directed unions. He has asked (again in \cite{Day1}) whether the 
class EG coincides with the class AG of all amenable groups 
(see also \cite{Cho}). 

   Today, we know that there are finitely generated groups in AG which are not
in EG; this has first been shown using growth estimates (\cite{Gri2},
\cite{Gri3}), and more recently by an elegant argument of Stepin 
(see \cite{Ste}, based on \cite{Gri2}). 
\par

   One knows also finitely {\it presented} groups in AG which are not in EG; 
more precisely, the finite presentation
$$
  G \, = \, 
  \left\langle \,\, a , b , c, d , t \,\,\, \bigg\vert \,\,\, 
\aligned
   & a^2 = b^2 = c^2 = d^2 = bcd = (ad)^4 = (adacac)^4 = 1  \,\, \\
   & t^{-1}at = aca 	\quad t^{-1}bt = d \quad t^{-1}ct = b 
          \quad t^{-1}dt = c \,\,
\endaligned
  \right\rangle
$$
defines an amenable group which is not elementary amenable 
(\cite{Gri6}, \cite{Gri7}).
\bigskip

\noindent  {\bf 12. Bartholdi's presentation.} It has later
been shown that the group $G$ of \cite{Gri6} has a presentation with two
generators only (namely $a$ and $t$) and four relations of total length
 $109 =  2 + 19 + 32 + 56$. Here are Bartholdi's computations, where $T$
stands for $t^{-1}$. \par

   The relations $c = aTata$, $d = tcT$ and $b = Tct$ show first that the
relations $c^2 = d^2 = b^2 = 1$ may be deleted in the presentation above,
and second that the generators $b,c,d$ may also be deleted. Thus
$$
  G \, = \, 
  \left\langle \,\, a ,  t \,\,\, \bigg\vert \,\,\, 
\aligned
    a^2 = TctctcT &= (atcT)^4 = (atcTacac)^4 = 1  \,\, \\
          T^2ct^2 &= tcT \,\,
\endaligned
  \right\rangle
$$
where $c$ holds for $aTata$. The relation $TctctcT = 1$ implies
$T^2ctctc = 1 = tcT^2ctc$ (by conjugation), hence also (using $c^{-1} = c$)
$$
   1 \, = \, \left( T^2ctctc \right)\left( tcT^2ctc\right)^{-1}
     \, = \, T^2ct^2 \left( tcT \right)^{-1}
$$
using free simplifications, so that the relation $T^2ct^2 = tcT$ may also be
deleted. Finally, one observes that $atcT$ is conjugate to
$Tatc = (Tata)^2$ so that
$(atcT)^4 = 1$ may be written $(Tata)^8 = 1$,
and one observes also that $atcTacac$ is equal to
$ataTataTaaTataaaTata$, so is conjugate to
$T^2ataTat^2aTata$. One obtains finally Bartholdi's presentation
$$
  G \, = \, 
  \left\langle \,\, a ,  t \,\,\, \bigg\vert \,\,\, 
    a^2 = TaTatataTatataTataT = (Tata)^8 = (T^2ataTat^2aTata)^4 = 1  \,\, 
  \right\rangle .
$$
\bigskip

\noindent {\bf 13. Categorical considerations.} 
For a given integer $k$, let $F_k$ be the free group on $k$ generators $\{ s_1 , \hdots , s_k \}$ and let
$X_k$ denote the space of all {\it marked groups on $k$ generators,} namely of 
all data $F_k \twoheadrightarrow \Gamma$, where $\twoheadrightarrow$ indicates a
homomorphism {\it onto.}  There is an appropriate topology on $X_k$, for which two
quotients  $\pi \colon F_k \twoheadrightarrow \Gamma$ and 
$\pi' \colon F_k \twoheadrightarrow \Gamma '$ are ``near''\,  
each other if the  corresponding Cayley graphs have balls of 
``large''\, radius around the unit element which are isomorphic.   
This makes $X_k$ a compact space; one shows for example that the closure of   
the subset of $X_k$ corresponding to finite groups contains the subset of $X_k$
corresponding to residually finite finitely generated groups. For details,  
see \cite{Gri2}, \cite{Cha} and \cite{Ste}. \smallskip 

   It would be interesting to find pairs $(Y,Z)$ where \par

   \qquad $\bullet$ $Y$ is a compact subspace of $X_k$, \par

   \qquad $\bullet$ $Z$ is a ``small''\, (e.g. countable) subset of 
$Y$, consisting of amenable groups, \par

   \qquad $\bullet$ $Y \setminus Z$ consists of non-elementary amenable groups,
or more generally \newline
      \phantom{ahahaha} the set of elementary amenable 
       groups in  $Y \setminus Z$ is  of first category.

\noindent The point is that the space $Y$ contains a dense $G_{\delta}$
consisting of amenable groups which are not elementary amenable.
(As usual a $G_{\delta}$ in $Y$ is a countable intersection of open subsets  
of $Y$.) \par

   One such pair has been constructed in \cite{Gri2} and analyzed in
\cite{Ste}, with $Z$ a countable set of virtually $2$-step solvable groups 
and with  $Y \setminus Z$  consisting of infinite torsion groups. 
Understanding other such pairs would probably help us understanding
the closures of $AG_k$ and of $EG_k$ in $X_k$, where
$AG_k$ [respectively $EG_k$] denotes the subspace of $X_k$ containing 
marked groups $\pi \colon F_k \twoheadrightarrow \Gamma$ with $\Gamma$ amenable
[resp. elementary amenable].
\bigskip

\noindent {\bf 14. Variation on one question of Day.}
 Let us denote by BG the smallest class of groups containing finitely
generated groups of subexponential growth (see Definition 64) and closed with
respect to the four operations of Day listed in 11 above, namely with
respect to  $(i)$ taking subgroups, $(ii)$ forming factor groups, $(iii)$ group
extensions  and $(iv)$ upwards directed unions. \par 
Question: {\it does one have BG=AG ?}
\bigskip

\noindent {\bf 15. Other definitions of amenability for groups;
topological groups.} The natural setting for amenability of groups
is that of topological groups, mainly locally compact groups. A substantial part
of the theory consists in showing the equivalence of a large number of
definitions. \par

   Let $G$ be a {\it Hausdorff topological} group. Denote by ${\mathcal C}^b (G)$
the Banach space of bounded  continuous functions from $G$ to $\C$, with the
supremum norm. For $\xi \in {\mathcal C}^b (G)$  and $g \in G$, let 
${}_g\xi \in {\mathcal C}^b (G)$ be the function  $x \to f(g^{-1}x)$. Denote by 
${\mathcal U}{\mathcal C}^b (G)$ the closed subspace of  ${\mathcal C}^b (G)$ of functions  
$\xi$ for which the mapping  $g \mapsto {}_g\xi$ from $G$ to ${\mathcal C}^b (G)$
is continuous. 
   The following are known to be equivalent 
(see Theorem 3 in \cite{Day2} and Theorem 4.2 in \cite{Ric2}): 
\begin{enumerate}[$\bullet$]
\item there exists a left-invariant mean on ${\mathcal U}{\mathcal C}^b (G)$,
\item any continuous action $G \times Q \to Q$ of $G$
by affine transformations of a non-empty compact convex subset $Q$ of a Hausdorff locally convex topological
vector space has a fixed point. 
\end{enumerate}

\noindent The group $G$ is {\it amenable} if these properties hold. In case
$G$ is assumed to be locally compact, here is a short list of other
equivalent properties: 
\begin{enumerate}[$\bullet$]
\item there exists a left-invariant mean on ${\mathcal C}^b (G)$,
\item there exists a left-invariant mean on $L^{\infty}(G)$,
\item the unit representation of $G$ is weakly contained in the left regular representation of $G$ on $L^2(G)$,
\item for any continuous action $G \times X \to X$ of $G$ by homeomorphisms of a non-empty compact space $X$, 
there exists a $G$-invariant probability measure on $X$.
\end{enumerate}

\noindent The last point, on $G$-invariant measures, goes back to a paper by
Bogolyubov, see \cite{Bogl}, quoted by Anosov \cite{Ano}. This paper,
published in Ukrainian in 1939, has remained unnoticed; the paper does {\it
not} quote von Neumann \cite{NeuJ}, and it is conceivable that Bogolyubov has
introduced independently the notion of amenability.  
About relations between amenability, growth and existence of invariant
measures, we would also like to quote \cite{Bekl}. \smallskip

   The list above is very far from being complete! (See 16; other
items could be: several formulations of the F\o lner property for locally compact
groups, the Reiter-Glicksberg property, the existence of approximate units in
the  Fourier algebra, $\ldots$) See, e.g., the books \cite{Gre1}, \cite{Pat} and
\cite{Wag},  as well as \cite[Chapter 8]{Rei}, \cite{Eym2}, \cite[Chapter 4]{Zim},  
\cite[in particular Theorem 10.11]{Wag} and  \cite[Chapter 2]{Lub}.
In case of a countable group (with the discrete topology), here is the most
recent characterization of amenability with which one of the authors has
been involved: a countable group $G$ is amenable if and only if, for any action
of $G$ by homeomorphisms  on the Cantor discontinuum $K$, there exists a
probability measure on $K$ which is invariant by $G$ \cite{GiH2}. \smallskip

We would like to point out that some attention has been given to
topological groups which are not locally compact
(in \cite[Section 4]{Ric2} among other places). 
For example, let ${\mathcal U}({\mathcal H})_{st}$ be the group of unitary operators 
on a separable infinite dimensional Hilbert space ${\mathcal H}$, with the {\it strong
topology;} then ${\mathcal U}({\mathcal H})_{st}$ is amenable, namely there exists a left
invariant mean  on ${\mathcal U}{\mathcal C}^b \left( {\mathcal U}({\mathcal H})_{st} \right)$, but
there does not exist any left invariant mean on 
${\mathcal C}^b \left( {\mathcal U}({\mathcal H})_{st} \right)$ \cite{Har1, Har2}.
Moreover, this group does have closed subgroups which are not amenable;
indeed, if  ${\mathcal H} = \ell^2(F_n)$ for a free group $F_n$ of rank  $n \ge 2$, 
then  ${\mathcal U}({\mathcal H})_{st}$ has clearly a discrete subgroup isomorphic to
$F_n$, as observed in \cite{Har3}.  
Here is another example involving non locally compact
topologies; let $G$ be the group of real points of an $\R$-algebraic group
and let $\Gamma$ be a subgroup of $G$ which is dense for the
{\it Zariski topology;} if $\Gamma$ is amenable, so is $G$
(see \cite{Moo}, and Theorem 4.1.15 in \cite{Zim}).
\medskip

   Let us  mention the following: for a locally compact group
$G$ which is almost connected (this means that the quotient of $G$ by the
connected component of $1$ is compact), the three properties 
\par\noindent
$\phantom{abcde}$   $G$ is amenable, \newline
$\phantom{abcde}$   $G$ does not contain a discrete subgroup which is free on 
                    $2$  generators, \newline
$\phantom{abcde}$  $G/r(G)$ is compact, \par\noindent 
are equivalent.
This is due to Rickert: Theorem 5.5 in \cite{Ric2}, building on \cite{Ric1};
see also Theorem 3.8 in \cite{Pat}. Recall that the {\it solvable radical}
$r(G)$ of a locally compact group $G$ is the largest connected closed normal
solvable subgroup of $G$ \cite{Iwa}. (One may define similarly the
{\it amenable radical} of $G$ as the largest amenable closed normal subgroup
of $G;$ see Lemma 1 of Section 4 in \cite{Day1} and
Proposition 4.1.12 in \cite{Zim}.)
\par 

   This result of Rickert reduces in some sense the problem of understanding the
class of  amenable locally compact groups to totally disconnected groups; we believe 
moreover that the  most important (and difficult) part of the problem is that
which concerns finitely generated groups.  
\bigskip

\noindent {\bf 16. Cohomological definitions of amenability.}
There are various (co)homological characterizations of amenability. \par

One is that of Johnson: a group $G$ is amenable if and only if
$H^1(\ell^1(G),M^*)$ is reduced to $\{0\}$ whenever $M^*$ is a $G$-module
dual to some Banach $G$-module $M$ \cite{Joh}. It follows that the bounded
cohomology of an amenable group is always reduced to $\{0\};$ this is given by
Gromov (Section 3.0 in \cite{Gro1}) together with a reference to an unpublished 
explanation of Philip Trauber - hence the name ``Trauber theorem''.
\par

   Another one is in terms of ``uniformly finite homology''; it
applies to finitely generated groups, and indeed to metric spaces in a much
broader class. Such a space $X$ is {\it not} amenable if and only if the group
$H_0^{uf}(X)$ is reduced to $\{0\}$ (in this statement, one may take $\R$ as
coefficients, or equivalently $\Z$); this is one way to express that the
F\o lner condition does {\it not} hold in $X$ \cite{BlW1}. \par

   It seems also appropriate to quote here a theorem of Brooks: let $G$
be the covering group of a normal covering $M$ of a compact manifold $X;$
then $G$ is amenable if and only if $0$ is in the spectrum of the
Laplace-Beltrami operator acting on the space of square-integrable functions
on $M$ (see \cite{Bro}, or the exposition in \cite{Lot}). \par

   There are other conditions in terms of other ``coarse''\,
(co)homology theories of the groups, or in terms of K-theory of appropriate 
algebras associated with the group (see various papers by G. Elek,
including \cite{Ele2}). \par

   Let us mention that there are interesting cohomological {\it consequences}
of amenability. For example, let $G$ be a group which has an Eilenberg-MacLane
space $K(G,1)$ which is a finite complex; if $G$ is amenable, then $G$ has
Euler characteristic $\chi (G) = 0$ (a particular case of Corollary 0.6
of Cheeger and Gromov \cite{ChGr}, who use $\ell^2$-cohomology methods, 
and also a result of B. Eckmann, who uses other methods \cite{Eck}). 
Also, let $G$ be the fundamental group of some closed $4$-manifold
$M;$ if $G$ is infinite and amenable, then $\chi(M) \ge 0$ \cite{Eck}.

\bigskip

\noindent {\bf 17. Variations on amenability of groups.}
   There are standard variations on the pseudogroup ${\mathcal G}_G$ and the notion
of amenability. 
\par

One is to consider the pseudogroup ${\mathcal G}_{G \times G}$
associated as in Example 2.(i) with the action of $G \times G$ on $G$ defined by
$(x,y) \circ g = xgy^{-1}$. It is classical that ${\mathcal G}_{G \times G}$
is amenable if and only if ${\mathcal G}_G$ is amenable. In other words: $G$ has a
left invariant mean if and only if $G$ has a two-sided invariant mean 
(Lemma 1.1.1 and Lemma 1.1.3 in \cite{Gre1}).
\par

Another variation is to consider the action of $G$ on $G \setminus \{1\}$
defined by $x\circ g = xgx^{-1}$ and the notion of {\it inner amenability} for
a group. It is obvious that an amenable group is inner amenable. 
Straightforward examples (such as non-trivial direct products of free groups
and amenable groups) show that there are non-amenable groups which are  
inner amenable. More on this in \cite{BeHa}, \cite{Eff}, \cite{GiH1} and
\cite{HS2}. \par

   A third variation is to consider a subgroup $H$ of $G$ and the pseudogroup
${\mathcal G}_{G/H}$ associated with the natural action of $G$ on $G/H$. The 
subgroup $H$ is said to be {\it co-amenable} in $G$ if ${\mathcal G}_{G/H}$ is
amenable. There is a comprehensive analysis of this notion in \cite{Eym1};  
see also \cite{Bekk}, in particular Theorem 2.3. In case $G = F_m$ is a free
group of finite rank, a criterion for co-amenability of a subgroup in terms of
{\it cogrowth} is given in \cite{Gri1} (see Item 52 below). 
One may generalize actions of $G$ on $G/H$ to actions of $G$ on locally
compact spaces; co-amenability of $H$ is then a particular case of a
notion of amenability for actions known as {\it amenability in the sense of
Greenleaf} \cite{Gre2}. \par

   The notion of amenability for a group and that of co-amenability for a 
subgroup may both be viewed as particular cases of a notion for $G$-mappings,
for which we refer to \cite{AnaR}. In case of a group $G$ with the discrete
topology, it can be defined as follows. Let $X,Y$ be two Borel spaces given
with measure classes $\mu,\nu$ and with actions of $G$  by non-singular 
invertible Borel mappings, and let $\phi \colon X \to Y$ be a surjective Borel 
mapping such that $\phi_*(\mu) = \nu;$ thus there is a canonical linear
isometric mapping by which we identify the Banach space $L^{\infty}(Y,\nu)$
with a closed  $G$-invariant subspace of $L^{\infty}(X,\mu)$. 
Say the mapping $\phi$ is {\it amenable} if
there exists a $G$-equivariant linear mapping  
$E \colon L^{\infty}(X,\mu) \to L^{\infty}(Y,\nu)$ 
which is a conditional expectation, namely which is positive
and which restricts to the identity on $L^{\infty}(Y,\nu)$. 
   Example~1: $X = G$ and $Y$ is reduced to one point; then $X \to Y$ is
amenable if and only if $G$ is amenable. Example 2: $X = G/H$ for a
subgroup $H$ of $G$ and $Y$ is reduced to a point; then $X \to Y$ is
amenable if and only if $H$ is co-amenable in $G$.
Example 3: $X = G \times Z$ for a $G$-space $Z$ (with $G$ acting from the
left on itself and diagonally on the product $G \times Z$); then the
projection  $G \times Z \to Z$ is amenable if and only if the action of $G$
on $Z$ is amenable {\it in the sense of Zimmer} \cite[Section 4.3]{Zim}.\par

   There are other notions, including the three following ones: 
{\it K-amenability} \cite{Cun}, {\it weak amenability} \`a la Cowling-Haagerup
\cite{CowH}, and {\it a-T-menability} \`a la Gromov. (See 7.A and 7.E in
\cite{Gro3}, and \cite{BekCV}; in fact Gromov rediscovered the class of groups
having ``Property 3B''\ of Akemann and Walter in \cite{AkWa}.)

\section*{II.4. Tarski number of paradoxical group actions}

\noindent
Consider more generally the pseudogroup ${\mathcal G}_{G,X}$ associated with a
group action $G \times X \to X$ (see again Example 2.(i)). 
\medskip

\noindent {\bf 18. Definition.} For $\gamma \colon S \to T$ in
${\mathcal G}_{G,X}$,  define the Tarski number of $\gamma$ as
the smallest  ``number of pieces''\, $n$ such that there exists a
partition  $S = \sqcup_{1 \le j \le n}S_j$ and elements $g_1 , \hdots , g_n$
in $G$ with 
$\gamma (x) = g_j(x)$ for all  $x \in S_j, j \in \{1 , \hdots , n\}$. \par

   The  {\it Tarski number} of a paradoxical  ${\mathcal G}_{G,X}$-decomposition
$$
   X \, = \, X_1 \bigsqcup X_2 \quad , \quad \gamma_1 \colon X_1 \to X
       \quad , \quad \gamma_2 \colon X_2 \to X
$$
as above is the sum of the Tarski number of $\gamma_1$ and of that of $\gamma_2$.
It is clear that such a sum is an integer $\ge 4$. \par

 When ${\mathcal G}_{G,X}$ is not amenable, we define the {\it Tarski number}
${\mathcal T}(G,X)$ of the action  $G \times X \to X$ as the minimum of the
Tarski numbers of the paradoxical ${\mathcal G}_{G,X}$-decompositions of $X;$ 
when ${\mathcal G}_{G,X}$ is amenable, we set ${\mathcal T}(G,X) = \infty$. 
For a group $G$ acting on itself by left multiplication, we write
${\mathcal T}(G)$ rather than ${\mathcal T}(G,G)$. 
\bigskip

\noindent {\bf 19. Observation.} Let $G$ be a group given
together with a subgroup $G'$ and a quotient group $G''$. 
It is straightforward that one has
$$
\aligned
   {\mathcal T}(G) \, &\le \, {\mathcal T}( G') \\
   {\mathcal T}(G) \, &\le \, {\mathcal T}(G'') .
\endaligned
$$
For example, for the first of these inequalities, view $G$ as a disjoint
union of cosets of $G'$. 
\par

Each group $G$ has a finitely generated subgroup $G'$
such that ${\mathcal T}(G') = {\mathcal T}(G)$.
Indeed, assuming $G$ to be non-amenable, consider a paradoxical decomposition
$$
  G \, = \, X_1 \sqcup \hdots \sqcup X_m \sqcup Y_1 \hdots \sqcup Y_n
     \, = \, g_1 X_1 \sqcup \hdots \sqcup g_mX_m 
    \, = \, h_1 Y_1 \sqcup \hdots \sqcup h_nY_n
$$
containing $m + n = {\mathcal T}(G)$ pieces
(where $X_1 , \hdots , X_m , Y_1 , \hdots , Y_n$ are subsets of $G$
and  $g_1 , \hdots , g_m , h_1 , \hdots , h_n$ are elements of $G$).
Let $G'$ be the subgroup of $G$ generated by 
$\{ g_1 , \hdots , g_m , h_1 , \hdots , h_n \}$.
Set $X'_i = X_i \cap G'$ for all $i \in \{1 , \hdots , m\}$ and
$Y'_j = Y_j \cap G'$ for all $j \in \{1 , \hdots , n\}$.
Then 
$$
  G' \, = \, X'_1 \sqcup \hdots \sqcup X'_m \sqcup Y'_1 \hdots \sqcup Y'_n
     \, = \, g_1 X'_1 \sqcup \hdots \sqcup g_m X'_m 
    \, = \, h_1 Y'_1 \sqcup \hdots \sqcup h_n Y'_n
$$
so that ${\mathcal T}(G') \le{\mathcal T}(G)$. With the first inequality of
the present observation, this shows that ${\mathcal T}(G') = {\mathcal T}(G)$.
(One may observe a fortiori that $X'_1 , \hdots , Y'_n$ are non-empty.)
It follows that one has
$$
   {\mathcal T}(G) = \inf \big({\mathcal T}(G') \big)
$$
where the infimum is taken over all finitely generated subgroups $G'$ of $G$.
\par

   It should be interesting to study how the Tarski number behaves with respect
to other group theoretical constructions such as extensions
and  HNN-constructions.  In particular, for the latter, we
have in mind some presentations of  the Richard Thompson's $F$
group \cite{CaFP}; recall that $F$ is a group which does not
have non-abelian free subgroups, which is an HNN-extension of itself \cite{BrGe},
that $F$ is inner-amenable \cite{Jol}, that $F$ has non-abelian free subsemigroups
so that it is not supramenable  (see Chapter V below), and
that one does not know whether $F$ is amenable or not. 
\bigskip

{\bf 20. Proposition (Jonsson, Dekker).} {\it For a group $G$,
one has ${\mathcal T}(G) = 4$ if and only if $G$ contains a non-abelian free
subgroup.}

\medskip

\noindent
{\it Proof.} For the free group $F_2$ on $2$ generators $g$ and $h$, it is
classical that ${\mathcal T}(F_2) = 4;$ 
see, e.g., Figure 4.1 in \cite{Wag}.
We recall this as follows. Set 
$$
\aligned
   A_1 \, &= \, W\big( g \big) \\
   A_2 \, &= \, W\big( g^{-1} \big) \\
   B_1 \, &= \, W\big( h \big) \cup \left\{ 1 , h^{-1} , h^{-2} , \hdots
             \right\} \\
   B_2 \, &= \, W\big( h^{-1} \big) \smallsetminus \left\{ h^{-1} , h^{-2} ,
              \hdots \right\} 
\endaligned
$$
where $W(x)$ denotes the subset of $F_2$ consisting of reduced words on
$\{g,h\}$ with $x$ as first letter on the left, for 
$x \in \{g,g^{-1},h,h^{-1}\}$. Then
$$
   F_2 \, = \, A_1 \bigsqcup A_2 \bigsqcup B_1 \bigsqcup B_2  \, = \, 
            A_1 \bigsqcup gA_2 \, = \, B_1 \bigsqcup hB_2 .
$$
It follows that ${\mathcal T}(F_2) = 4$. \par
   Observation 19 shows that
${\mathcal T}(G) = 4$ for any group $G$ containing a subgroup isomorphic to $F_2$.
\smallskip

   Conversely, let $G$ be a group with ${\mathcal T}(G) = 4$, so
that there exist subsets $X_1 , X_2 , Y_1 , Y_2$ and elements
$g_1 , g_2 , h_1 , h_2$ in $G$ such that
$$
  G \, = \, X_1 \bigsqcup X_2 \bigsqcup Y_1 \bigsqcup Y_2   \, = \, 
      g_1X_1 \bigsqcup g_2X_2 \, = \, h_1Y_1 \bigsqcup h_2Y_2 .
$$
Set $g = g_1^{-1} g_2$ and $h = h_1^{-1} h_2$.
Then, one has successively
$$
\aligned
   X_1 \, &= \, G \smallsetminus gX_2 \, = \, gX_1 \bigsqcup  gY_1 \bigsqcup gY_2 \\
   X_1 \, &\supset \, gX_1 \, \supset \, \hdots \, \supset \, g^{k-1}X_1 \,
     \supset  g^kY_j \qquad (k \ge 1 \quad \text{and} \quad j = 1,2) \\
   X_2 \, &= \, G \smallsetminus g^{-1}X_1 \, = \, g^{-1}X_2  \bigsqcup g^{-1}Y_1
       \bigsqcup g^{-1}Y_2 \\
   X_2 \, &\supset \, g^{-1}X_2 \, \supset \, \hdots \, \supset
       g^{-k+1}X_2 \, \supset \,  g^{-k}Y_j 
      \qquad (k \ge 1 \quad \text{and} \quad j = 1,2)
\endaligned 
$$
so that
$$
   g^kY_j \, \subset X_1 \cup X_2 \qquad \text{for all} \qquad
        k \in \Z \, , \, k \ne 0 \quad \text{and} \quad j = 1,2 .
$$
One has similarly
$$
   h^kX_j \, \subset Y_1 \cup Y_2 \qquad \text{for all} \qquad
        k \in \Z \, , \, k \ne 0 \quad \text{and} \quad j = 1,2 .
$$
Hence $g$ and $h$ generate in $G$ a free subgroup of rank $2$, by a
classical lemma going back essentially to F. Klein, and sometimes known as
the ``table-tennis lemma''\, (see, e.g., \cite{Har4}). \par
   The argument above is our rephrasing of the proof of Theorem 4.8 in
\cite{Wag}. \hfill \qed 

\bigskip

Proposition 20 is an unpublished work from the 40's by B. Jonsson (a
student of Tarski) and is a particular case of results of Dekker published in
the 50's. For precise references, see the Notes of Chapter 4 in \cite{Wag}.
\par
   Let us also mention that, for a group $G$ containing a non abelian free group
and for an action $G \times X \to X$ with stabilizers 
$\{g \in G \, \vert \, gx = x \}$ which are abelian for all $x \in X$, 
the corresponding Tarski number is also given by 
${\mathcal T}(G,X) = 4$ (Theorem 4.5 in \cite{Wag}).\bigskip

\noindent 
{\bf 21. Proposition.} {\it For a non-amenable torsion group 
$G$, one has ${\mathcal T}(G) \ge 6$.}

\medskip

\noindent
{\it Proof.} By Proposition 20 we know that ${\mathcal T}(G) \ge 5$. We assume
that ${\mathcal T}(G) = 5$, and we will reach a contradiction. \par
The hypothesis implies that there exist 
subsets $X_1 , X_2 , Y_1 , Y_2 , Y_3$ and elements
$g_1 , g_2 , h_1 , h_2 , h_3$ in $G$ such that
$$
  G \, = \, X_1 \bigsqcup X_2 \bigsqcup Y_1 \bigsqcup Y_2  \bigsqcup Y_3 \, = \, 
    g_1X_1 \bigsqcup g_2X_2 \, = \, h_1Y_1 \bigsqcup h_2Y_2 \bigsqcup h_3Y_3.
$$
Let $n$ denote the order of $g \Doteq g_1^{-1}g_2$. As in the proof of
Proposition 11, one has 
$$
 X_1 \, \supset \, gX_1 \, \supset \, \hdots \, \supset \, g^{n-1}X_1 \, 
   \supset \, g^n \left( Y_1 \bigsqcup Y_2 \bigsqcup Y_3 \right) .
$$
But now $g^n = 1$ and this is absurd. Hence ${\mathcal T}(G) > 5$. \hfill \qed
\bigskip

\noindent {\bf 22. Question.}  Does there exist a 
group $G$ with Tarski number ${\mathcal T}(G)$ equal to 5 ? to 6 ? More generally, what
are the possible values of ${\mathcal T}(G)$ ?

\section*{II.5. F\o lner condition for pseudogroups}

\noindent
Let $({\mathcal G} , X)$ be a pseudogroup of transformations. For a subset
${\mathcal R}$ of ${\mathcal G}$ and a subset $A$ of $X$, we define the
{\it ${\mathcal R}$-boundary of $A$} as
$$
\partial_{\mathcal R}A \, = \, \left\{ \, x \in X \smallsetminus A \, \Bigg\vert \,
\aligned
   &\text{there exists $\rho \in {\mathcal R}\cup {\mathcal R}^{-1}$ such that} \\
   &\qquad  x \in \alpha(\rho) \quad \text{ and } \quad \rho(x) \in A \, 
\endaligned \right\} .
$$

\noindent
{\bf 23. Definition.} {\it The pseudogroup $({\mathcal G} , X)$  satisfies the 
{\it F\o lner condition} if 
$$
\aligned
   &\text{for any finite subset ${\mathcal R}$ of ${\mathcal G}$ and for any real number
            $\epsilon > 0$} \\
   &\text{there exists a finite non-empty subset $F = F({\mathcal R}, \epsilon)$
            of $X$} \\
   &\text{such that} \quad \vert \partial_{\mathcal R} F \vert \, < \,
                \epsilon \vert F \vert
\endaligned
$$
where $\vert F \vert$ denotes the cardinality of the set $F$.}
\bigskip

\noindent {\bf 24. Ahlfors and F\o lner.}
    Ideas underlying the F\o lner condition go back at least to Ahlfors.
(F\o lner does {\it not} refer to this work.) Ahlfors defines an open
Riemann surface $S$ to be {\it regularly exhaustible} if, for some
appropriate metric $g$ in the conformal class defined by the complex
structure of $S$, there exists a nested sequence 
$\Omega_1 \subset \Omega_2 \subset \hdots$ of domains with smooth boundaries
such that $\bigcup_{n \ge 1} \Omega_n$ is the whole surface and such that
$$
 \lim_{n \to \infty} \frac{
   \vert \partial \Omega_n \vert_g}{ \vert \Omega_n \vert_g} \, = \, 0 
$$
where $\vert \Omega \vert_g$ denotes the area of a domain $\Omega$ 
and where $\vert \partial \Omega \vert_g$ denotes the length of its boundary,
both with respect to $g$. (A lemma of Ahlfors shows that this does not depend
on the choice of $g$.) These sequences may be used to define averaging processes,
as Ahlfors did first and as F\o lner did later. \par

    Using this notion, Ahlfors has developped a geometric approach to the
Nevanlinna theory of distribution of values of meromorphic functions, known as
Ahlfors theory of covering surfaces. In particular, he gave a generalization of
the second main theorem of Nevanlinna on defect. 
(See Section 25 in Chapter III of \cite{Ahl}; see also Chapter XIII in
\cite{Nev}, Chapter 5 in \cite{Hay}, Theorem 6.5 on page 1223 of
\cite{Oss},  \cite{Sto} and \cite{ZoKe}.)
\par

Amenability of coverings of Riemann surfaces can also be expressed in terms
of Teichm\"uller spaces \cite{McM2}. \par
\bigskip

\noindent{\bf 25. Theorem.} {\it A pseudogroup of transformations is amenable if and
only if it satisfes the F\o lner condition.}\bigskip

F\o lner's original proof (for a group acting on itself by left
multiplications) goes back to 1955 \cite{Fol}. The proof has been simplified
by Namioka \cite{Nam} (who generalized F\o lner's result to one-sided
cancellative semigroups), and extended to group actions by Rosenblatt
\cite{Ros1}; the best place to read it is probably Section 2.1 of \cite{Co1}. 
In case of a group $G$ acting by conjugation on $G \smallsetminus \{1\}$, the
proof can also be found in \cite{BeHa}, and it applies {\it verbatim} to an
action of $G$ on any set $X$. All these references use essentially techniques
of functional analysis. (See also Wagon's comment about the implication 
$(6) \Longrightarrow (1)$ in Theorem 10.11 of \cite{Wag}.) \par

The proof below, in Items 26 and 36, uses completely different techniques. 

\bigskip

\noindent {\bf 26. Beginning of the proof of Theorem 25.} We prove here   
the implication ``F\o lner condition $\Rightarrow$ existence of an
invariant mean''. \par

   Let ${\mathcal M}(X)$ denote the set of all means on $X$, namely of all 
finitely additive probability measures on $X$ (see Conditions (fa) and (no)
in Definition 4). Let $\ell^{\infty} (X)$ denote the Banach
space of all bounded functions on $X$, with the norm of uniform convergence;
it is standard\footnote{See
footnote $37$ in \cite{NeuJ}, where von Neumann refers in turn
to Lebesgue's ``Le\c cons sur l'int\'egration''\, 
(1905).
}
that ${\mathcal M}(X)$ can be identified with a subset of the unit
ball in the dual space of $\ell^{\infty} (X)$. It is also standard that the
weak$^*$-topology makes ${\mathcal M}(X)$ into a compact space. \par

For each finite non-empty subset $F$ of $X$, we consider the mean
$$
\mu_F \colon \quad \left\{ \aligned
   {\mathcal P}(X) \qquad &\longrightarrow \qquad\quad [0,1] \\
    A \quad\qquad &\longmapsto \qquad 
           \frac{ \vert A \cap F \vert }{ \vert F \vert} 
                      \endaligned \right.
$$
in ${\mathcal M}(X)$. Consider also the set
$$
\mathcal N \, = \, \left\{ 
\hskip.2cm (\mathcal R, \epsilon) \hskip.2cm 
\vert
\mathcal R \subset \mathcal G \hskip.2cm \text{is finite, and} \hskip.2cm
\epsilon \in \R , \hskip.1cm \epsilon > 0
\right\}
$$
ordered by
$$
   ({\mathcal R}, \epsilon) \, \le \, ({\mathcal R}' , \epsilon')
   \qquad \text{if} \qquad {\mathcal R}\subset {\mathcal R}' \quad \text{and} \quad
   \epsilon \ge \epsilon '.
$$
Notation being as in Definition 23 of the F\o lner condition (which is now
assumed to hold), 
\[
\begin{split}
 (*) & \ \ \ \ \ \ \ \ \ \ \ \ \ \ \ \ \ \ \ \ \ \ \ \ \ \ \ \ \ \ \ \ \ \ \ \  \left( \mu_{F({\mathcal R}, \epsilon)} \right)_{ ({\mathcal R}, \epsilon) \in {\mathcal N}} \ \ \ \ \ \ \ \ \ \ \ \ \ \ \ \ \ \ \ \ \ \ \ \ \ \ \ \ \ \ \ \ \ \ \ \ \ \ \ \ \ \ \ \ \ \ \ \ \ \ \ \ \ \ \ \ \ \ \ \ \ \ \ \ \ \ \ \ \ \ \ \ \ \ \ \ \ \ \ \ \ \ \ \ \ \ \ \ \ \ \ \ \ \ \ \ \ \ \ \ \ \ \ \ \ \ \ \
\end{split}
\]
becomes a {\it net.} By compacity of ${\mathcal M}(X)$, this net has
a cluster point, say $\mu$ (we use the terminology of \cite[Chapter 2]{Kel}).
The proof consists in showing that $\mu$ is ${\mathcal G}$-invariant; in other
words, given a subset $A$ of $X$ and a transformation $\gamma$ in ${\mathcal G}$
with $A \subset \alpha(\gamma)$, one has to show that 
$\mu\big( \gamma (A) \big) = \mu (A)$. \par

   We choose a number $\delta > 0$. As $\mu$ is a cluster point of the family
$(*)$, there exists $({\mathcal R}, \epsilon) \in {\mathcal N}$ such that
$$
\aligned
   (i) \qquad &({\mathcal R}, \epsilon) \, \ge \, (\{ \gamma \},\delta),
        \quad i.e., \quad {\mathcal R}\ni \gamma \quad \text{and} \quad 
        \epsilon \le \delta , \\
   (ii) \qquad &\vert \mu_{F({\mathcal R}, \epsilon)}(A) - \mu(A) \vert 
         \, \le \, \delta , \\
   (iii) \qquad &\big\vert \mu_{F({\mathcal R}, \epsilon)}\big( \gamma (A) \big)
          - \mu \big( \gamma (A) \big) \big\vert \, \le \, \delta.
\endaligned
$$
From now on, we write $F$ instead of $F({\mathcal R}, \epsilon)$. Define
$$
\aligned
   A_{i,i} \, &= \, \left\{ \, a \in A \, \vert \, a \in F 
           \quad \text{and} \quad \gamma(a) \in F \, \right\}  =  A \cap F \cap \gamma^{-1}(F), \\
   A_{i,o} \, &= \, \left\{ \, a \in A \, \vert \, a \in F 
           \quad \text{and} \quad \gamma(a) \in \partial_{{\mathcal R}}F 
                     \, \right\} = A \cap F \cap \gamma^{-1}(X \smallsetminus F), \\
   A_{o,i} \, &= \, \left\{ \, a \in A \, \vert \, a \in \partial_{{\mathcal R}}F 
           \quad \text{and} \quad \gamma(a) \in F
                     \, \right\} = A \cap (X \smallsetminus F) \cap \gamma^{-1}(F), \\
   A_{o,o} \, &= \, \left\{ \, a \in A \, \vert \, a \notin F 
           \quad \text{and} \quad \gamma(a) \notin F
                     \, \right\} = A \cap (X \smallsetminus F) \cap \gamma^{-1}(X \smallsetminus F)
\endaligned
$$
(think of ``inside''\, for ``$i$''\, and of ``outside''\, for ``$o$'').
Observe that 
$A = A_{i,i} \sqcup A_{i,o} \sqcup A_{o,i} \sqcup A_{o,o}$,
with the first three sets being finite. Observe also that
$$
\aligned
   (iv) \qquad &A \cap F \, = \, 
          A_{i,i} \sqcup 	A_{i,o} \qquad \text{so that}\qquad
          \vert A \cap F \vert \, = \, 
          \vert A_{i,i} \vert + \vert A_{i,o} \vert \\
   (v) \qquad &\gamma \quad \text{induces a bijection} \quad
           A_{i,i} \sqcup A_{o,i} \to \gamma (A) \cap F \\
   &\text{so that} \qquad 
            \vert \gamma(A) \cap F \vert \, = \, 
            \vert A_{i,i} \vert + \vert A_{o,i} \vert \\
   (vi) \qquad &\partial_{\mathcal R}F \, \supset \, 
                \partial_{\{ \gamma,\gamma^{-1} \}}F \, \supset \,
                \gamma(A_{i,o}) \cup A_{o,i} \\
   &\text{so that} \qquad
             \vert A_{i,o} \vert + \vert A_{o,i} \vert \, \le \,
             2 \vert \partial_{\mathcal R}F \vert \, \le \, 
             2 \epsilon \vert F \vert .
\endaligned
$$
It follows from $(iv)$ to $(vi)$ that
$$
   (vii) \qquad
   \big\vert \vert \gamma (A) \cap F \vert - \vert A \cap F \vert \big\vert
   \, \le \, 2 \epsilon \vert F \vert .
$$
Using the definition of the mean $\mu_F$ and the conclusion of the F\o lner
condition, one may rewrite $(vii)$ as
$$
   (viii) \qquad
   \big\vert \mu_F \big( \gamma (A) \big) - \mu_F(A) \big\vert \, \le \, 
   2 \epsilon
$$
so that one obtains finally
$$
 \big\vert \mu \big( \gamma (A) \big) - \mu(A) \big\vert \, \le \,
   2 \delta + 2 \epsilon \, \le \, 4 \delta
$$
using $(ii)$, $(iii)$ and $(viii)$. As the choice of $\delta$ is arbitrary,
this ends the proof of one implication of Theorem 25. \hfill \qed

\bigskip

 \noindent {\bf 27. Remark.}  In case of a locally finite graph $X$ with
finitely many orbits of vertices under the full automorphism group 
(for example in case of a Cayley graph), F\o lner condition 
is equivalent to the existence of a nested sequence 
$F_1 \subset F_2 \subset \hdots$ of finite subsets of the vertex set $X^0$
such that $\cup_{n \ge 1}F_n = X^0$ and 
$\lim_{n \to \infty} \vert \partial F_n \vert / \vert F_n \vert = 0;$ 
see our Section III.2 for amenable graphs and for
the notation $\partial F_n$, and Theorem 4.39 in
\cite{Soa} for the equivalence. \par

   In the case of a group $G$ acting on a set $X$, the
F\o lner condition is most often expressed in a way involving the symmetric
difference between a finite subset $F$ of $X$ and its image $gF$ by some 
$g \in G;$ for the equivalence of this with the analogue of our Definition
23, see Proposition 4.3 in \cite{Ros1}. \par

   For groups, F\o lner condition implies the existence of F\o lner sets with
extra tiling properties, and this is useful for showing extensions to amenable
groups of the Rohlin theorem from ergodic theory \cite{OrWe}.
   
\section*{{\bf III. Amenability and paradoxical decompositions for metric spaces}}

\section*{III.1. Gromov condition and doubling condition}

\noindent
Let $X$ be a metric space and let $d$ denote the distance on $X$. \par

 For $S,T \subset X$, a mapping 
$\phi \colon S \to T$ ({\it not} necessarily a bijection) is a {\it bounded
perturbation of the identity} if  $\sup_{x \in S} d(\phi(x) , x) < \infty$.
We will denote by 
$$ 
            {\mathcal B}(X)
$$ 
the collection of all these maps.
(This would be an example of a ``pseudo-semigroup'', but we will
not use this term again below.) \par

   As in Example 2.(iii), we denote by ${\mathcal W}(X)$ the pseudogroup
of all bijections, between subsets of $X$, which are bounded perturbations of
the identity. \par

   For a subset $A$ of $X$ and a real number $k > 0$, we denote by
$$
   {\mathcal N}_k(A) \, = \, \{ \, x \in X \, \vert \, d(x,A) \le k \, \}
$$
the {\it $k$-neighbourhood} of $A$ in $X$. 
\par

   Recall that a metric space is \emph{locally finite}\footnote{The 
terminology ``discrete'' of the 1999 published version
is not appropriate. More on this in Chapter VI.} 
if its subsets of finite diameter are finite. \bigskip

\noindent {\bf 28. Definitions.}
    A locally finite metric space $X$ is said to be  {\it amenable} 
[respectively {\it paradoxical}] if the pseudogroup ${\mathcal W}(X)$ is amenable
[resp. paradoxical].\bigskip

   {\it Caution.} This definition is not convenient for non-locally finite metric
spaces, because the pseudogroup ${\mathcal W}(\R)$ is paradoxical.
Indeed, the bijections
$$
\gamma_{even} \, \colon \, \bigcup_{n \in \Z}[2n,2n+1[ \, \longrightarrow \, 
            \R 
\qquad \text{and} \qquad 
\gamma_{odd} \, \colon \,  \bigcup_{n \in \Z}[2n+1,2n+2[ \, \longrightarrow \, 
            \R 
$$
defined by $\gamma_{even}\vert_{[2n, 2n+1[}(t) = 2t - (2n + 1/2)$
and $\gamma_{odd}\vert_{[2n+1, 2n+2[}(t) = 2t - (2n + 3/2)$,
are in ${\mathcal W}(\R)$ and define a paradoxical decomposition of $\R$.
\par

   A notion of amenability for {\it some} non-locally finite metric spaces is
suggested in Remark 42. \bigskip

\noindent {\bf 29. Definition.} 
A locally finite metric space $X$ is said to satisfy the {\it Gromov condition} if 
there exists a mapping $\phi \colon X \to X$ in ${\mathcal B}(X)$  
such that 
$$
          \big\vert \phi ^{-1} (x) \big\vert \ge 2
$$ 
for all $x \in X$. 
\medskip

   This terminology refers in particular to the ``lemme 6.17''\,
in \cite{GrLP}, introduced there as ``le meilleur moyen de montrer qu'un
groupe est non-moyennable''; see also Item 0.5$.C''_1$ in \cite{Gro3}.
\bigskip

\noindent {\bf 30. Definition.} The locally finite metric space $X$ satisfies the
{\it doubling condition} if there exists a constant $K > 0$ such that
$$
   \big\vert {\mathcal N}_K(F) \big\vert \, \ge \, 2 \vert F \vert
$$
for any non-empty finite subset $F$ of $X$.
\par
   It is of course equivalent to ask that there exists a constant $k > 0$
and a number $\epsilon > 0$ such that 
$$
   \big\vert {\mathcal N}_k(F) \big\vert \, \ge \, (1 + \epsilon) \vert F \vert
$$
for any non-empty finite subset $F$ of $X;$ indeed, this implies
$\vert {\mathcal N}_K (F) \vert \ge 2 \vert F \vert$ for any non-empty
finite subset $F$ of $X$, with $K = nk$ and $n$ an integer such that
$(1 + \epsilon)^n \ge 2$.
\bigskip

   \noindent {\bf 31. Bipartite graphs and matchings.} 
Let $B = Bip(Y,Z;E)$ be a {\it bipartite} graph with two classes $Y,Z$ of
vertices and with edge set $E;$ by definition of ``bipartite'', any
edge $e \in E$ is incident with one vertex in $Y$ and one vertex in $Z;$ we
consider here simple graphs, namely graphs without loops and {\it without
multiple edges.} Recall that, for integers $k,l \ge 1$, a {\it perfect
$(k,l)$-matching} of $B$ is a subset $M$ of $E$ such that any $y \in Y$ 
[respectively any $z \in Z$] is incident to exactly $k$ edges in $M$ 
[resp. $l$ edges in $M$]. \par

   For a set $F$ of vertices of $B$, we denote by $\partial_E F$ the set of
vertices in $B$ which are not in $F$, and are connected to some vertex of $F$
by some $e \in E$. 

   Let again $X$ be a metric space, as earlier in the present section. With two
subsets $S,T \subset X$ and a real number $K \ge 0$, one associates the
bipartite graph $B_{K}(S,T)$ with vertex classes $S$ and $T$, and with an edge
connecting $x \in S$ and $y \in T$ whenever $d(x,y) \le K$; note that, by definition,
$S$ and $T$ are disjoint in the vertex set of $B_K(S, T)$, even if they need not be as subsets of $X$. 
 Observe that $X$ is locally finite if and only if $B_K(X,X)$ is locally finite for all $K \ge 0$.
\bigskip

 {\bf 32. Theorem.} {\it For a locally finite metric space $X$, the following
conditions are equivalent (with ${\mathcal B}(X)$ as before Definition 28). 
\smallskip

   $(i)$ The space $X$ is paradoxical. \par

   $(ii)$ There exists a mapping $\phi \colon X \to X$ in ${\mathcal B}(X)$  
such that $\big\vert \phi^{-1} (x) \big\vert = 2$ for all $x \in X$. \par

   $(iii)$ There exists a mapping $\phi \colon X \to X$ in ${\mathcal B}(X)$  
such that $\big\vert \phi^{-1} (x) \big\vert \ge 2$ for all $x \in X$ 
(namely $X$ satisfies the Gromov condition). \par

   $(iv)$ The space $X$ satisfies the doubling condition. \par

   $(v)$ There exists a real number $K > 0$ for which the bipartite graph
$B_K(X,X)$ has a perfect $(2,1)$-matching. \par

   $(vi)$ The pseudogroup ${\mathcal W}(X)$ does not satisfy the F\o lner condition.}
	\bigskip

\noindent {\bf 33. Observations.} As there are amenable groups of exponential
growth, for example finitely generated solvable groups which are not 
virtually nilpotent, Conditions $(ii)$ and $(iii)$ are not connected to
growth, as suggested in \cite{DeSS}, but indeed to amenability, as already
observed in our Introduction. \par
   For a recent survey on growth and related matters, see \cite{GriH}. \par
   Some of the implications of Theorem 32 may be made more precise. See for
example Proposition 54 below. 
\bigskip

\noindent  {\bf 34. Proof of Theorem 32.} \par

   $(i) \Longleftrightarrow (ii)$. If $X$ is paradoxical, there exists a
partition $X = X_1 \sqcup X_2$ and two bijections 
$\gamma_j \colon X_j \to X$ in ${\mathcal W}(X)$. The mapping $\phi \colon X \to X$ defined
by $\phi (x) = \gamma_j(x)$ for $x \in X_j \, (j=1,2)$ satisfies $(ii)$.
\par

   Conversely, given a mapping $\phi \colon X \to X$ as in $(ii)$, one uses the
axiom of choice to order the two points of $\phi^{-1}(x)$ for each 
$x \in X$, say as $\phi^{-1}(x) = \left( \gamma_1^{-1}(x) , \gamma_2^{-1}(x) \right)$.
This provides a paradoxical decomposition involving the mappings $\gamma_1$ and
$\gamma_2$. 

\medskip

    The implications $(ii) \Longrightarrow (iii) \Longrightarrow (iv)$
are straightforward. 
   Condition $(v)$ is nothing but a rephrasing of Condition $(ii)$. 
	
\medskip	

   $(vi) \Longrightarrow (iv)$.  If ${\mathcal W}(X)$ does not satisfy the
F\o lner condition, there exists $\epsilon > 0$ and a non-empty finite subset
${\mathcal R}$ of ${\mathcal W}(X)$ such that, for any non-empty finite subset $F$ of
$X$, one has 
$\vert F \cup \partial_{\mathcal R} F \vert \ge (1 + \epsilon) \vert F \vert $. 
Setting 
$$
    C \, = \, \max_{\rho \in {\mathcal R}\cup {\mathcal R}^{-1}} 
       \sup_{x \in \alpha (\rho)} d(\rho(x),x)
$$
(see Definition 1 for the notation $\alpha(\rho)$), one has a fortiori
$$
   \big\vert {\mathcal N}_C(F) \big \vert \, \ge \, (1 + \epsilon) \vert F \vert 
$$
for any non-empty finite subset $F$ of $X$.

\medskip

   $(i) \Longrightarrow (vi)$. The contraposition
$\text{not}(vi) \Longrightarrow \text{not}(i)$ may be checked as
follows: if the pseudogroup ${\mathcal W}(X)$ 
satisfies the F\o lner condition, it is amenable by Proof 26, 
so that ${\mathcal W}(X)$ is not paradoxical by the straightforward part of
the Tarski alternative (Remark 6.(i)). 

\medskip

   We have now shown all but the right lowest $\Rightarrow$ in the 
following diagram: 
$$
\begin{array}{*{20}c}
 & & & & (v) &&&& (vi) &&
 \\
 &&&&&&&&&&
 \\
&&&& \Updownarrow &&&& \Downarrow &&
 \\
&&&&&&&&&&
 \\
 (vi) & \Leftarrow & (i) & \Leftrightarrow & (ii) & \Rightarrow & (iii) & \Rightarrow & (iv)
 & \Rightarrow & (v)
\end{array}
$$
   
   For the last implication $(iv) \Longrightarrow (v)$,
we follow \cite{DeSS} and call upon a form of the Hall-Rado Theorem.
More precisely, with the notation of Theorem 35 below and with $k = K$,
$(iv)$ implies that  
$\big\vert \partial_E F \big\vert  \ge  2 \vert F \vert$
for any subset $F$ of $Y$ or of $Z$, so that $(v)$ follows. \hfill \qed
\bigskip

   All what we will need about the Hall-Rado theorem can be found in
\cite{Mir} but, as a first background, we recommend also the discussion in 
Section III.3 of \cite{Bol}. (Recall that ``Hall''\,\, refers to
{\it Philip} Hall.)
\bigskip

\noindent
{\bf 35 Theorem (Hall-Rado).} {\it Let $B = Bip(Y,Z;E)$ be a locally finite
bipartite graph and let $k \ge 1$ be an integer. Assume that one has
$$
\aligned
   &\big\vert \partial_E F \big\vert \, \ge \, k \vert F \vert
      \qquad \text{for all finite subsets $F$ of $Y$} \\
   &\big\vert \partial_E F \big\vert \, \ge \, \phantom{k } \vert F \vert
      \qquad \text{for all finite subsets $F$ of $Z$.}
\endaligned
$$
Then there exists a perfect $(k,1)$-matching of $B$. }
\medskip

\noindent {\it On the proof.} Consider the bipartite graph 
$B_k = B\left( \sqcup_{1 \le j \le k}Y_j , Z ; E_k \right)$ where 
$\sqcup_{1 \le j \le k}Y_j$ denotes a disjoint union of $k$ copies of
$Y$, and where, for each edge $e \in E$ with ends $y \in Y$ and  
$z \in Z$, there is one edge $e_j \in E_k$ with ends the vertex $y_j \in Y_j$
corresponding to $y$ and the vertex $z$, this for each $j \in \{1 \hdots k\}$.

One the one hand, the hypothesis implies that
$$
   \big\vert \partial_{E_k} F \big\vert \, \ge \,  \vert F \vert
$$
for all finite subset $F$ of $\sqcup_{1 \le j \le k}Y_j$ or of $Z$.
On the other hand, 
there exists a perfect $(k,1)$-matching of $B$ if and only if
there exists a perfect $(1,1)$-matching of $B_k$. 
It follows that one may assume $k = 1$ without loss of generality. 
\smallskip

   By the most usual form of the Hall-Rado theorem, 
there are subsets
$M_Y , M_Z$ of $E$ such that the edges in $M_Y$ [respectively in $M_Z$] are
pairwise disjoint, and such that each $y \in Y$ [resp. each $z \in Z$] is
incident with exactly one edge in $M_Y$ [resp. in $M_Z$]; see, e.g., Theorem
4.2.1 in \cite{Mir}. Thus $M_Y \cup M_Z$ define a spanning subgraph of $B$
whose connected components are either edges, or simple polygons with a number
of edges which is even and at least $4$, or infinite lines. 
(This argument is standard: see e.g. the middle of page 317 in \cite{Nas}.)
\par

   One may color the edges of the latter subgraph in black and white such 
that each vertex of $B$ is incident to exactly one black edge. The set of
black edges thus obtained is a perfect $(1,1)$-matching of $B$. \hfill \qed
\bigskip

   If $k = 1$, observe that the condition of the Theorem is also necessary for
the existence of a perfect $(1,1)$-matching. If $k \ge 2$, it is not so (consider
a complete bipartite graph with $\vert Y \vert = 1$ and $\vert Z \vert = k$),
despite the statement following Definition 6 of \cite{DeSS}.
\bigskip

\noindent {\bf 36. End of proof of Theorem 25.} We show here the implication
``existence of an invariant mean $\Rightarrow$ F\o lner condition'', or rather its contraposition: 
we assume that $( {\mathcal G} , X )$ does not satisfy the F\o lner condition,
and we have to prove that $X$ has no ${\mathcal G}$-invariant mean. \smallskip

   {\it First case: } $X$ is a metric space and ${\mathcal G}$ is the pseudogroup
${\mathcal W}(X)$. 
   Implication $(vi) \Longrightarrow (i)$ of Theorem 32 shows that
$X$ is paradoxical, hence that $X$ is not amenable. The
proof of Theorem 25 is complete in this case. \smallskip

   {\it General case.} If $( {\mathcal G} , X )$ does not satisfy the F\o lner
condition, there exists a number $\epsilon > 0$ and a non-empty finite subset 
${\mathcal R}$ of ${\mathcal G}$ such that
$$
   \vert \partial _{\mathcal R} F \vert \, > \, \epsilon \vert F \vert
$$
for any non-empty finite subset $F$ of $X$. Define a metric $d_{\mathcal R}$ on
$X$ by
$$
   d_{\mathcal R}(x,y) \, = \, \min \left\{ \, n \in \N \, \Bigg\vert
\aligned
   &\text{there exists} \quad \rho_1 , \hdots , \rho_n \in {\mathcal R}\cup
      {\mathcal R}^{-1} \quad \text{such that} \\
   &\rho_n \Big( \rho_{n-1} \big ( \hdots \rho_1 (x) \hdots \big) \Big)
      \, \text{is defined and is equal to} \, y
\endaligned 
\, \right\} 
$$
with the understanding that $d_{\mathcal R}(x,y) = \infty$ if there exists no
such $n$. One has a posteriori
$$
  \vert {\mathcal N}_1 (F) \vert \, \ge \, (1 + \epsilon) \vert F \vert
$$
for any non-empty finite subset $F$ of $X$, where the neighborhood 
${\mathcal N}_1 (F)$ refers to the metric $d_{\mathcal R}$ (for the definition
of ${\mathcal N}_1$, see before Definition 28). Hence the pseudogroup 
${\mathcal W}(X, d_{\mathcal R})$ is not amenable by the previous case. As
${\mathcal W}(X, d_{\mathcal R}) \subset {\mathcal G}$, the pseudogroup ${\mathcal G}$ itself is
not amenable either. \hfill \qed

\bigskip

\noindent {\bf 37. Definition.} Recall that two metric spaces $X,Y$ are {\it
quasi-isometric} if there exist constants $\lambda \ge 1 \, , \, C \ge 0$ and a
mapping $\phi \colon X \to Y$ such that
$$
   \frac{1}{\lambda} d (x_1 , x_2)  - C \, \le \, 
   d\big( \phi(x_1) , \phi(x_2) \big) \, \le \, 
   \lambda d (x_1 , x_2)  + C 
$$
for all $x_1 , x_2 \in X$ and
$$
   d \big( y , \phi(X) \big) \, \le \, C
$$
for all $y \in Y$. \par
   Recall also that $X$ and $Y$ are {\it Lipschitz equivalent} if there
exists a constant $\lambda \ge 1$ and a {\it bijection} $\psi \colon X \to Y$ such
that
$$
   \frac{1}{\lambda} d (x_1 , x_2)   \, \le \, 
   d\big( \psi(x_1) , \psi(x_2) \big) \, \le \, 
   \lambda d (x_1 , x_2)  
$$
for all $x_1 , x_2 \in X$. (See also Item 0.2.C in \cite{Gro3}.)
\bigskip

\noindent \textbf{38.\ Proposition.}
\emph{Let $X$ and $Y$ be two uniformly locally finite metric spaces
which are quasi-isometric. Then $X$ is amenable [respectively paradoxical]
if and only if $Y$ is so.}

\vskip.2cm

\noindent \emph{Proof.}
For uniformly locally finite\footnote{We are grateful to Volker Diekert for having pointed out to us
the omission of uniform local finiteness in the hypotheses in our previous version.}
metric spaces, the Gromov condition of Definition 29
is clearly invariant by quasi-isometry. \hfill \qed

\bigskip

\noindent {\bf 39. Examples.} For each prime $p$, there are uncountably many
$2$-generated $p$-groups which are amenable and pairwise {\it not}
quasi-isometric; see  \cite{Gri2} for $p = 2$
and \cite{Gri3} for $p \ge 2$. \bigskip

\noindent {\bf 40. Examples.} There are uncountably many
$2$-generated torsion-free groups which are paradoxical and pairwise {\it not}
quasi-isometric \cite{Bow}.
\bigskip  

\noindent {\bf 41. Remark.}    
It is a result due independently to Volodymyr Nekrashevych \cite{Nek1} 
and Kevin Whyte \cite{Why} that two uniformly discrete non-amenable metric spaces $X$ and 
$Y$  are quasi-isometric if and only if they are Lipschitz equivalent. This  
answers a question of Gromov (Item $1.A'$ in \cite{Gro3}); see also \cite{Pap}    
and \cite{Bogp} for partial answers. 
\bigskip

\noindent {\bf 42. Remark.} Let $(\Omega , d_{\Omega})$ be a metric space. 
A subset $X$ of $\Omega$ is a {\it separated net} if 
there exists a constant $r > 0$ for which the
two following properties hold:    
(i) $d_{\Omega}(x,y) \ge r$ for all $x,y \in X$, $x \ne y$, and 
(ii) $X$ is a  maximal subset of $\Omega$  for this property
(this implies $d_{\Omega}(\omega,X) \le 2r$ for all 
$\omega \in \Omega$). Such nets exist by Zorn's Lemma. \par

   If the metric space $\Omega$ is ``slim and well-behaved''\, in the
sense of \cite{MaMT}, for example if $\Omega$ is a Riemannian manifold with
Ricci curvature bounded from below and the injectivity radius of the exponential 
map positive, then two nets in $\Omega$ are quasi-isometric to each other. (See
Theorems 3.3 and 3.4 in \cite{MaMT}, as well as \cite{Kan1}, \cite{Kan2} and
\cite{Nek1}, \cite{Nek2}.)
   For such slim and well-behaved spaces, there are natural notions of
amenability  and paradoxes, defined via their nets; this has appeared in several
places, including \cite{BlW1}.
   Proposition 38 carries over to these spaces, by definition.  
\bigskip

\noindent {\bf 43. Examples. } There are uncountably many Riemann
surfaces of constant curvature $-1$ which are amenable as metric
spaces, and which are pairwise {\it not} quasi-isometric \cite{Gri5}.
\bigskip

\section*{III.2. Graphs as metric spaces, isoperimetric constants}

\noindent
Let $X = (X^0,X^1)$ be a graph with vertex set $X^0$ and with edge set
$X^1$ (say $X$ has no loops and no multiple edges, for simplicity). 
If $X$ is connected, $X^0$ is naturally a uniformly discrete metric space, the distance $d(x,y)$ between 
two vertices  $x,y \in X^0$ being the minimal number of edges in a path between 
them. \par

   For a disconnected graph $X$, there are also notions of combinatorial
distances. For example, if $X$ is a subgraph of a  connected graph $Y$ which is
clear from the context, one may restrict to $X^0$ the distance defined on $Y^0$ as
above. One may also set $d(x,y) = \infty$ for $x,y$ in different connected
components of $X$. 

In this section, we assume that $X = (X^0, X^1)$ is a graph given together with
a metric $d \colon X^0 \times X^0 \longrightarrow {\mathbf R}_+$ such that
$d(x,y)$ is the combinatorial distance whenever $x,y$ are vertices
in the same connected component of $X$.

\bigskip

\noindent {\bf 44. Definition.}
A locally finite graph $X$ is said to be {\it amenable} or {\it
paradoxical} if the metric space $X^0$ is so in the sense of Definition 28.
\bigskip

   For a subset $F$ of $X^0$, the  boundary $\partial_E F$
defined in graph theoretical terms in Item 31 (here $E = X^1$)
coincides with ${\mathcal N}_1 (F) \setminus F$, where ${\mathcal N}_1 (F)$ is the
neighborhood defined in metrical terms before Definition 28. We will write
$$
   \partial F \, = \, {\mathcal N}_1 (F) \setminus F
$$
below.
\bigskip 

\noindent {\bf 45. Definition.} The 
{\it isoperimetric constant} of the graph $X$ is
$$
   \iota (X) \, = \, \inf \left\{ \, 
 \frac{ \vert \partial F \vertÊ}{ \vert F \vert } \, \Big\vert \, 
   F \subset X^0 \quad \text{is finite and non-empty } \, \right\} .
$$
For example, 
$\iota(X) = 0$ as soon as $X$ is finite.
\bigskip

\noindent {\bf 46. Variations.} There are several variations on the
definition of the isoperimetric constant in the literature, because a boundary
$\partial F$  could be defined using \smallskip

   either  vertices {\it outside} $F$ as here (before Definition 45)
          or in \cite{BeSc} and \cite{McM1}, \par
   or vertices {\it inside} $F$ as in \cite{Dod} or \cite{CoSa}, \par
   or vertices {\it both} inside and outside $F$ as in \cite[page 24]{OrWe}, 
\par
   or  {\it edges} connecting vertices inside $F$ to those outside $F$ 
as in \cite{BiMS} or \cite{Kai1}. 
\smallskip

   For example, denoting by $\partial _* F$ the set of {\it edges}
connecting a vertex of $F$ to a vertex outside $F$, there is another
isoperimetric constant
$$
   \iota_* (X) \, = \, \inf \left\{ \, 
 \frac{ \vert \partial_* F \vertÊ}{ \vert F \vert } \, \Big\vert \, 
   F \subset X^0 \quad \text{is finite and non-empty } \, \right\} .
$$
for the graph $X$. One has $\iota_*(X) \ge \iota(X);$ if $X$ has
maximal degree $k$, one has also $\iota_*(X) \le k \iota(X)$.

\bigskip

\noindent{\bf 47. Example} {\it Let $d$ be an integer, $d \ge 3$.
For a tree $T$ in which every vertex is of degree at least $d$,
the isoperimetric constant satisfies the inequality
$$
      \iota (T) \, \ge \, d-2.
$$
If $T$ is regular of degree $d$, then $\iota (T) = d-2$.} 
\medskip

\noindent{\it Proof.} As we have not found a convenient published reference for this
very standard fact, we indicate now a proof. We denote by $T(d)$ the regular tree of degree $d$.
\par

Let $F$ be a finite subset of the vertex set of $T$, let $X$ denote
the subgraph of $T$ induced by $F$, let $X_1 , \hdots , X_N$ denote its
connected components, and let $F_i$ denote the vertex set of $X_i$, for  
$i \in \{1 , \hdots , N \}$. We claim that
$$
    \vert \partial F \vert  \, \ge \,  (d-2)\vert F \vert + 2 .
$$\par

   Assume first that $X$ is connected. We proceed by induction on 
$\vert F \vert$. If $\vert F \vert = 1$, then 
$\vert \partial F \vert \ge d = (d-2)\vert F \vert + 2$ and the claim is
obvious. Assume now that  $\vert F \vert = k \ge 2;$ let $y \in F$ be a
vertex of $X$-degree $1$, and let $Y$ be the subgraph of $X$ induced by 
$F \setminus \{y\}$. One has
$$
   \vert \partial F \vert \, \ge \, 
   \vert \partial ( F \setminus \{y\} ) \vert + d-2 \, \overset{*}{\ge} \, 
   (d-2) \big( \vert F \vert - 1 \big) + 2 + d-2 \, = \, 
   (d-2) \vert F \vert + 2
$$
where $\overset{*}{\ge}$ holds because of the induction hypothesis.
(It is easy to check that 
$\vert \partial F \vert = (d-2) \vert F \vert + 2$ in case $T = T(d)$.) 
\par

   Assume now that $X$ has $N \ge 2$  connected components, and proceed by
induction on $N$. As $T$ is a tree, one may assume the enumeration of the
$F_i$ 's such that $\partial F_1$ has at most one vertex in common with
$\partial \left( \bigcup_{2 \le i \le N}  F_i \right)$. Then
$$
   \big\vert \partial F \big\vert
   \, \ge \,
   \vert \partial F_1 \vert \, + \, 
   \big\vert \partial \Big( \bigcup_{2 \le i \le N} F_i \Big) \big\vert
   \, - \, 1
   \, \overset{**}{\ge} \, 
   (d-2) \vert F_1 \vert + 2 + (d-2) \sum_{i=2}^N \vert F_i \vert + 2 - 1
   \, > \, 
   (d-2) \vert F \vert + 2
$$
where $\overset{**}{\ge}$ holds because of the induction hypothesis.
\par

   It follows that $\iota (T) \ge d-2$, with equality for a $d$-regular tree. \hfill \qed
	
\bigskip

   Recall that a {\it hanging chain} of length $k$ in a graph $X$ is a path 
of length $k$ (with $k+1$ vertices, $k-1$ so-called inner ones and the two 
end-vertices) with all inner vertices of degree $2$ in $X$. It is obvious
that, if $X$ has hanging chains of arbitrarily large lengths, then 
$\iota (X) = 0$. The following is a kind of converse, for trees.
\bigskip

   \noindent{\bf 48. Example} {\it Let $T$ be a connected infinite locally finite 
tree without end-vertices and let $k$ be an integer, $k \ge 2$. 

If $T$ has no hanging chain of length $> k$, then
$$
   \iota (T) \, \ge \, \frac{1}{2k} .
$$
Also $\iota (T) = 0$ if and only if T has arbitrarily long hanging chains.}
\medskip

\noindent {\it Proof:} see the proof of Corollary 4.2 in \cite{DeSS}. \hfill \qed

\bigskip

   Other interesting estimates of isoperimetric constants appear, for example,
in Section 4 of \cite{McM1}.
\bigskip

   \noindent {\bf 49. Definitions.} On a locally finite graph $X$, there is a
natural {\it simple random walk} with corresponding {\it Markov operator}
$T$. Suppose for simplicity that $X$ is connected and of bounded degree.
Consider the Hilbert space $\ell ^2 (X^0, \, deg )$ of functions $h$ from $X^0$
to $\C$ such that $\sum_{x \in X^0} deg(x) \vert h(x) \vert ^2 < \infty$, 
and the bounded self-adjoint operator $T$ defined on this Hilbert space by 
$$
   (Th)(x) \, = \, \frac {1}{deg(x)} \sum_{y \sim x} h(y)
$$
for $h \in \ell^2(X^0, \, deg )$, $x \in X^0$,  where $y \sim x$ indicates a
summation over the neighbours $y$ of the vertex $x$. The {\it spectral radius}
of $X$ is
$$
\aligned
   \rho(X) \, &= \, \sup \left\{ \, \langle h \vert Th \rangle \, 
        \big\vert \, h \in \ell ^2 (X^0, \, deg ) \, , \, \|h\|_2 \le 1 \, \right\} \\
   &= \,\sup \left\{ \, \vert \lambda \vert \,\,\, \big\vert \quad  \lambda
         \quad  \text{is in the spectrum of} \quad T \, \right\} .
\endaligned
$$
Observe that $1 - T$ is a natural analogue on $X$ of a Laplacian, so that
$1 - \rho (X)$ is often referred to as the first eigenvalue of the
Laplacian or (more appropriately) as the bottom of its spectrum. \par

   It is also known that, for a real number $\lambda$, the following are
equivalent:
$$
\aligned
(i)\phantom{iiiiii}
   &\text{there exists} \quad F \colon X^0 \to ]0,\infty[ \quad \text{such that}
         \quad \frac{1}{deg (x)} \sum_{y \sim x} F(y) = \lambda F(x) ,\\
(ii)\phantom{iiiii}
   &\text{there exists} \quad F \colon X^0 \to ]0,\infty[ \quad \text{such that}
         \quad \frac{1}{deg (x)} \sum_{y \sim x} F(y) \le \lambda F(x) ,\\
(iii)\phantom{iiii}
   &\text{one has} \quad \lambda \ge \rho (X),
\endaligned
$$
so that $(i)$ and $(ii)$ indicate alternative definitions
of the spectral radius.
In terms of the Laplace operator, $(i)$ and $(ii)$ are respectively
conditions about $(1-\lambda)$-harmonic and $(1-\lambda)$-superharmonic
functions.
(For a proof in terms of graphs, see Proposition 1.5 in \cite{DoKa}.  
But there are earlier proofs in the literature on irreducible stationary 
discrete Markov chains. The equivalence of (ii) and (iii) is standard;  the
equivalence with (i) is more delicate: \cite{Harr} and \cite{Pru}.) 
\par

   For $x,y \in X^0$ and for an integer $n \ge 0$, denote by $p^{(n)}(x,y)$ the
probability that a simple random walk starting at $x$ is at $y$ after $n$ 
steps. Then one has also 
$$\rho(X) = \limsup_{n \to \infty}\root n \of { p^{(n)}(x,y) };$$
in particular, the value of this $\limsup$ is independent on $x$ and $y$.
From this probabilistic interpretation of $\rho(X)$, one deduces easily that,
for a connected graph $X$ which is regular of degree $d \ge 2$,
one has $\rho(X) \ge 2 \sqrt{d-1} / d;$
equality holds if and only if $X$ is a tree.
\smallskip

(More generally, for any transition kernel  $p \colon X^0 \times X^0 \to [0,\infty[$ 
with reversible measure $\mu \colon X^0 \to ]0,\infty[$, so that $\sum_{z \in X^0}p(x,z) = 1$ and
$\mu(x) p(x,y) = p(y,x)\mu(y)$ for all $x,y \in X^0$, 
one introduces the Hilbert space $\ell^2(X^0,\mu)$,
and the self-adjoint operator $T$ defined by the kernel $p$ on 
$\ell ^2(X^0,\mu)$. Then the norm of $T$ is again equal to
$\limsup_{n \to \infty}\root n \of { p^{(n)}(x,y) }$.)
 \bigskip

\noindent{\bf 50. Lemma (an isoperimetric inequality).} {\it For a graph
$X$ which is regular of degree $d \ge 2$, one has 
$$
    \iota (X) \, \ge \, 4 \, \frac{ 1 - \rho(X)  } { \rho(X) }.
$$
}
\medskip

\noindent{\it Proof.} 
    Let ${\mathbb X} ^1$ denote the set of oriented edges of $X$.
(If $X$ is finite, the cardinal of ${\mathbb X}^1$ is twice the number of
geometric edges of $X$.) Each $e \in {\mathbb X}^1$ has a head
$e_+ \in X^0$ and a tail $e_- \in X^0$. For a function 
$h \in \ell^2 (X^0 , deg)$ with real values, one has
$$
   \left\langle h \vert Th \right\rangle \, = \, 
   \sum_{x \in X^0} h(x) \sum_{y \sim x} h(y) \, = \, 
   \sum_{ e \in {\mathbb X} ^1} h(e_+)h(e_-) \, = \, 
   \|h\|^2 - \frac{1}{2} \sum_{e \in {\mathbb X}^1} 
                 \Big( h(e_+) - h(e_-) \Big)^2 .
$$
  Let now $F$ be a finite non-empty subset of $X^0$, with boundary 
$\partial F$. Consider the function $h \in \ell^2(X^0,deg)$
defined by
$$
   h(x) \, = \, \left\{ 
\aligned
   \frac{1}{\sqrt{d}} \, \quad 
                &\text{if} \quad x \in F  \\
   \frac{1}{2\sqrt{d}} \quad 
                &\text{if} \quad x \in \partial F \\
   0 \,\,\, \quad &\text{otherwise}
\endaligned \right.
$$
One has clearly
$$
 (*) \ \ \ \ \ \ \ \ \  \|h\|^2 \, = \, \vert F \vert + \frac{1}{4} \vert \partial F \vert
          \, \ge \, 
         \vert F \vert \left( 1 + \frac{ \iota(X)}{4} \right) .  
$$
One has also
$$
    \frac{1}{2} \sum_{e \in {\mathbb X}^1}  \Big( h(e_+) - h(e_-) \Big)^2
   \, = \, \sum_{y \in \partial F} \,  
            \sum_{x \sim y} \Big( h(y) - h(x) \Big) ^2 
   \, \le \,  \vert \partial F \vert \, d \, \frac{1}{4d} .
$$
Together with $(^*)$, this implies that
$$
   \rho(X) \, \ge \, \frac{\left\langle h \vert Th \right\rangle}{\|h\|^2}
   \, \ge \, 1 \, - \, \frac{\vert \partial F \vert}
       {4 \vert F \vert \left( 1 + \frac{ \iota(X)}{4} \right)} .
$$
Taking the infimum over $\frac{ \vert \partial F \vert }{\vert F \vert}$
one obtains
$$
     \rho(X) \, \ge \, 1 \, - \, 
     \frac{ \frac{\iota(X)}{4}}{1  +  \frac{\iota(X)}{4} }
$$
and the lemma follows. \hfill \qed
\bigskip

   The previous lemma appears in several places (see N$^o$ 51 below). It is
related to Theorem 3.1 of \cite{BiMS}, which is stated in terms of the
constant 
$\iota_* (X)$ of our Item 46, and which shows that
$\iota_* (X) \ge 4(1 - \rho(X))$. Recently, T. Smirnova-Nagnibeda 
has improved the latter to
$$
   \iota_* (X) \, \ge \, \frac {d^2}{d-1} ( 1 - \rho(X) )
$$
(the improvement comes from choosing a test-function, playing the role
of the function $h$ in the proof above, which is more efficient than the one
chosen in \cite{BiMS}).
\smallskip

   For a {\it majoration} of $\iota(X)$ in terms of $1 - \rho(X)$ and $d$
(namely for an analogue of the ``Cheeger's inequality''),
see Theorem 2.3 in \cite{Dod} or Theorem 3.2 in \cite{BiMS}
(in each case with normalizations different from ours).
\bigskip

   \noindent{\bf 51. Theorem.} {\it Let $X$ be a connected graph which is of
bounded degree. The following are equivalent:
$$
\begin{aligned}
 &(i) \qquad    && X \quad \text{is paradoxical} \qquad  && \text{(see Definition 44),} \\
 &(ii) \qquad   && \iota(X) > 0   \qquad       &&\text{(see Definition 45),} \\
 &(iii)  \qquad && \rho(X) < 1    \qquad      & &\text{(see Definition 49),} \\
 &(iv) \qquad   && p^{(n)}(x,y) \, = \, o(\sigma ^n) \qquad &&  \text{for some } \sigma \in ]0,1[ \text{ and for all }   x,y \in X^0 
\end{aligned}
$$
\noindent and they imply that
$$
\begin{aligned} 
& (v) \qquad & \text{ the simple random walk on $X$ is transient.} \qquad && \qquad { \ }
\end{aligned}
$$ }
\medskip

\noindent {\it On the proof.} The equivalence $(i) \Longleftrightarrow (ii)$ is a
reformulation of Theorem 25  on the F\o lner condition. \par

   The equivalence $(ii) \Longleftrightarrow (iii)$ may be viewed as a discrete
analogue of the Cheeger-Buser inequalities for Riemannian manifolds
\cite{Che}, \cite{Bus}. For graphs as in the present theorem, it can be found in
\cite{Dod}, \cite{Var}, \cite{DoKe}, \cite{DoKa}, \cite{Ger}, \cite{Anc},
\cite{Kai1}; there are also similar arguments showing appropriate estimates for
{\it finite} graphs in several papers by Alon et alii, quoted
in \cite{Lub} (in particular near Propositions 4.2.4 and 4.2.5). \par

   For $(iii) \Longleftrightarrow (iv)$ and for equivalence with other 
conditions, see Theorem 4.27 in Soardi's notes on Networks \cite{Soa}.

   The implication $(iii) \Longrightarrow (v)$ is obvious.

      For groups, the equivalence 
$$
\text{amenability} \qquad \Longleftrightarrow \qquad \rho(X) = 1
$$
goes back to the pioneering papers of Kesten \cite{Kes1}, \cite{Kes2}. See
also \cite{Day3} and the review in \cite{Woe}. \hfill \qed
\bigskip

There are other conditions equivalent to $(i)$ to $(iv)$ above, for example
in terms of norms of Markov operators on $\ell^p$-spaces; see
\cite{Kai1}. \smallskip

   For locally finite graphs which are not necessarily of bounded degree, 
one has to modify some of the definitions above. Thus, for a finite set
$F$ of vertices of a graph $X$, one considers the sum $\|F\|$ of the degrees
of the vertices in $F$, the number $\|{\partial F}\|$ of edges with one end in
$F$ and the other end outside $F$, and the infimum $\tilde \iota (X)$ of the
quotients $\|{\partial F}\| / \|F\|$ (compare with Definition 45). 
For graphs of bounded degree, one has 
$\tilde \iota (X) = 0 \Longleftrightarrow \iota(X) = 0$, but in general\footnote{
Here is an example shown to us by Vadim Kaimanovich. Let 
$\left( h_j \right)_{j \ge 1}$ be a sequence of integers, all at least $2$, and
consider first a rooted tree $Y$ in which a vertex at distance $n$ of the root
is of degree
$$
\left\{
\aligned 
     \quad  k+2\phantom{3} \quad &\text{if} \quad n = \sum_{j=1}^k h_j \quad 
                 \text{for some} \quad k \ge 1, \\
     \quad  3\phantom{k+2} \quad   &\text{otherwise.}
\endaligned \right.
$$
Consider then the graph $X$ obtained from $Y$ by adding, for each vertex $x$ of
$Y$ at distance $n = \sum_{j=1}^k h_j$ from the root (for some $k$), the 
$\frac{1}{2}(k+1)(k+2)$ edges between the successors of $x$ in $Y$.
Then one has $\iota (X) > 0$ (because $Y$ is a spanning tree for $X$) and
$\tilde \iota (X) = 0$ (because $X$ contains induced subgraphs which are
complete graphs on $k+2$ vertices for $k$ arbitrarily large).
One has also $\rho(X) = 1$.
} on may
have $\tilde \iota (X) = 0$ and $\iota (X) > 0$. 
By a particular case of a result of Kaimanovich (Theorem 5.1 in \cite{Kai1}), 
one has $\tilde \iota (X) > 0 \Longleftrightarrow \rho (X) < 1$.
\par

   Graphs of unbounded degree are also covered by the arguments in \cite{DoKa}
and \cite{DoKe}. \bigskip

   Graphs give rise to several kinds of algebras, and it is a natural question in
each case to ask how the properties of Theorem 51 translate. For Gromov's 
{\it translation algebras} (see the end of 8.C$_2$ in \cite{Gro3}), there is a
hint in \cite{Ele1}. For other algebras associated with graphs (and more 
generally with oriented graphs), see \cite{KPRR} and \cite{KPR}.  Amenable
properties of certain kind of graphs (more precisely of bipartite graphs with
appropriate weights) are also important in the study of subfactors; see various
works by S. Popa, including \cite{Pop1} and \cite{Pop2}.\par
   Amenability has of course been one of the most important notions in the theory of
operator algebras since the works of von Neumann. We will not discuss more of
this here, but only refer to \cite{Co2} and \cite{Hel}.

\section*{{\bf IV. Estimates of Tarski numbers}}

\section*{IV.1. From relative growth to Tarski number of paradoxical decompositions}

 \noindent  
Let $G$ be a finitely generated group, given as a quotient
$$
    \pi \, \colon \, F_m \longrightarrow G
$$
of the free group $F_m$ on $m$ generators $s_1 , \hdots , s_m$,
for some $m \ge 1$. The purpose of the present section is to review 
notions which will be used in IV.2.
\bigskip

\noindent {\bf 52. Recall: relative growth, spectral radius and isoperimetric constant.} 
Let $\ell \colon F_m \to \N$
denote the word length on $F_m$ with respect to $s_1 , \hdots, s_m$. For each
integer $k \ge 0$, let $\sigma(\ker(\pi),k)$ denote the cardinality of the set  
$\left\{ \, w \in \ker(\pi) \, \vert \, \ell (w) = k \, \right\} $.
The {\it relative growth} of $\ker(\pi)$ 
(some authors say ``the cogrowth of $G$''!) is,
by definition,
$$
 \alpha_{\ker(\pi)} \, = \, \limsup_{k \to \infty} \root k \of { \sigma(\ker(\pi),k) } .
$$

If $\ker (\pi) \ne \{1\}$, it is easy to check that $\sqrt{2m-1} \le \alpha_{\ker (\pi)} \le 2m-1$,
and one shows more precisely that $\sqrt{2m-1} < \alpha_{\ker (\pi)}$,
see  \cite{Gri1}.
\smallskip

   The corresponding Cayley graph (with vertex set $G$ and with an edge between
two vertices $x,y$ if and only if $\ell(xy^{-1}) = 1$) has a spectral radius given
by the formula
$$
\rho \, = \, \left\{
\aligned
   \frac{\sqrt{2m-1}}{m} \hskip2.5cm
  \hskip.5cm  &\text{if} \hskip.5cm
   1 \le \alpha \le \sqrt{2m-1} 
\\
   \frac{\sqrt{2m-1}}{2m} \left( \frac{\sqrt{2m-1}}{\alpha} \, + \, 
               \frac{\alpha}{\sqrt{2m-1}} \right) 
   \hskip.5cm &\text{if} \hskip.5cm
     \sqrt{2m-1} < \alpha \le 2m-1 
\endaligned \right.
$$
\cite{Gri1}. It follows that the three conditions
$$
\aligned
   &\alpha \, = \, 2m - 1 \\
   &\rho \, = \, 1 \\
   &G \quad \text{is amenable}
\endaligned
$$
are equivalent; the equivalence of the last two is due to Kesten,
as already recalled in the proof of Theorem 51. 
   (In the present setting for the formula giving $\rho$ as a function of
$\alpha$, one has  $1 \le \alpha \le \sqrt{2m-1}$ if and only if $\alpha = 1$, 
if and only if  $\ker(\pi) = \{1\};$ but the formula makes sense and is correct 
for subgroups of $F_m$ which need not be normal, and then the range
$1 \le \alpha \le \sqrt{2m-1}$ is meaningful.) 
\bigskip

\noindent {\bf 53. Isoperimetric constant and doubling characteristic distance.}   
Let $X$ be a graph, with its set $X^0$ of vertices viewed as a metric space for 
the combinatorial distance $d$ as in Section III.2. 
A {\it doubling characteristic distance} for $X$ is (if it exists) an integer $K$
for which the doubling condition of Definition 30  holds, namely an integer $K$ 
such that  $$
   \vert {\mathcal N}_K (F) \vert \, \ge \, 2 \vert F \vert
$$
for any non-empty finite subset $F$ of $X^0$.
If the isoperimetric constant $\iota(X)$ of Definition 45 is strictly positive, 
the integer
$$
   K_X \,  = \, \left\lceil \frac{ \log 2 }{Ê\log ( 1 + \iota(X)) } \right\rceil
$$
is clearly a doubling characteristic distance, 
where $\lceil t \rceil $ indicates the least integer larger than or equal to
$t$.
\bigskip

{\bf 54. Proposition.} {\it Let $X$ be a graph with isoperimetric constant
$\iota(X) > 0;$ define $K_X$ as in the previous paragraph. Then there exists a
paradoxical decomposition involving a partition $X^0 = X^0_1 \sqcup X^0_2$ 
and two bounded perturbations of the identity $\phi_i \colon X^0_j \to X^0$
in ${\mathcal W}(X^0)$ such that
$$
   \sup_{x \in X^0_j} d(\phi_j(x) , x) \, \le \, K_X \qquad (j = 1,2).
$$}
\medskip

\noindent{\it Proof:} this is a quantitative phrasing of the
implication $(iv) \Longrightarrow (i)$ of Theorem 32, and follows from our Proof 
34. \hfill \qed
\bigskip

\noindent {\bf 55. Four functions.} Let $m$ be an integer, $m \ge 2$. \par

   For $\alpha \in ]\sqrt{2m-1},2m-1]$, set 
       $\rho_m(\alpha) \, = \, 
       \frac{\sqrt{2m-1}}{2m} \left( \frac{\sqrt{2m-1}}{\alpha} \, + \, 
               \frac{\alpha}{\sqrt{2m-1}} \right) \in 
               \left] \frac{\sqrt{2m-1}}{m},1\right]$. \par

   For $\rho \in ]0,1]$, set 
       $\iota(\rho) = 4 \frac{1 - \rhoÊ}{\rho} \in [0,\infty[$. \par

   For $\iota \in [0,\infty[$, set 
       $K(\iota) = \left\lceil \frac{ \log 2 }{Ê\log ( 1 + \iota) } \right\rceil
            \in \{1 , 2 , 3 , \hdots , \infty\}$,  
       with $\left\lceil \hdots \right\rceil$ as in 53. \par

   For $K \in \{1 , 2 , 3 , \hdots , \infty\}$, set
       $b_m(K) = \frac{m(2m-1)^K - 1}{m-1} $. \par

\noindent Observe that $\alpha \mapsto \rho_m(\alpha)$ and $K \mapsto b_m(K)$ are
increasing, while $\rho \mapsto \iota(\rho)$ and $\iota \mapsto K(\iota)$ are
decreasing. Observe also that, in the Cayley graph  of a group $G$ with respect to a
set of $m$ generators, a ball of radius $K$ has at most $b_m(K)$ elements, and
precisely $b_m(K)$ elements in case $G$ is free on $m$ generators.

\bigskip

\noindent{\bf 56. Theorem.} {\it Let $G = F_m/N$ be a group given as a quotient of the free
group on $m$ generators by a normal subgroup $N \ne \{1\}$. Let $\alpha_G$ denote 
the corresponding relative growth and let $\iota(X)$ denote the isoperimetric constant  
of the corresponding Cayley graph $X$ (see Definition 45 and Item 52). Using the
notation of the previous number, one has: \par

   (i) if $\alpha_G \le \alpha$ for some $\alpha \le 2m-1$, the Tarski number of $G$
satisfies
$$
  {\mathcal T}(G) \, \le \, 
       2 b_m \bigg( K \Big( \iota \left( \rho_m(\alpha) \right) \Big) \bigg) ,
$$

   (ii) if $\iota(X) \ge \iota$ for some $\iota \ge 0$, then
$$
  {\mathcal T}(G) \, \le \, 2 b_m \big( K (\iota) \big) .
$$}
\medskip

\noindent{\it Proof.} For (i), one has $\iota(X) \ge \iota\left( \rho_m(\alpha) \right)$
by the formula of Item 52 and by the isoperimetric inequality of Lemma 50,
and this implies $K_X \le K\left( \iota\left(\rho_m(\alpha)\right) \right)$.
If $\phi_j \colon X^0_j \to X^0$ are as in Proposition 54, one may write $X^0_j$
as a finite disjoint union of the sets
$$
   A_{j,g} \, = \, \left\{ x \in X^0_j \mid \phi_j(x) = gx \right\}
$$
for $g$ in the ball $B^G(K_X) = \left\{ g \in G \mid \ell(g) \le K_X \right\}$
(compare with Observation 9), this for $j = 1$ and $j = 2$.
As $\vert B^G(K_X) \vert \le b_m(K_X)$,
this ends the proof of (i). The end of the argument shows also (ii). \hfill \qed

\bigskip

\noindent {\bf 57. Comments and examples.}
Observe that we have argued with the Cayley graph of $G$ related to the
{\it right}-invariant distance $d(x,y) = \ell\left(xy^{-1}\right)$ on $G$, so that
the {\it left}-multiplications $x \mapsto gx$ are bounded perturbations of the
identity. \smallskip

   Let us now test the inequalities of Theorem 56. \smallskip

   (i) Let $F_2$ denote the free group of rank $2$ and let $X$ denote the Cayley
graph of $F_2$ with respect to some free basis ($X$ is of course a regular tree of
degree $4$). Kesten \cite{Kes1} has computed the spectral value of the corresponding
simple random walk as $ \rho(X)  =  \frac{ \sqrt 3 }{2}   \approx  0.86603 $ 
so that
$ \iota (X)   \ge  4 \frac{ 1 - \rho(X) }{\rho(X)}   \approx   0.6188$.
Hence
$ K_X   =  \left\lceil \frac{ \log 2 }{Ê\log ( 1.6188 )} \right\rceil =  2$ 
is a doubling characteristic distance. The resulting estimate
$$
{\mathcal T}(F_2) \, \le \, 2 \vert B^{F_2}(2) \vert 
    \, = \, 2 \big( 2^.3^2 - 1 \big)
    \, = \, 34
$$
compares rather badly with the correct value ${\mathcal T}(F_2) = 4$.\par

   A similar computation with the Cayley graph $Y$ of $F_3$ with respect to a free
basis gives $\rho (Y) = \frac{ \sqrt{5}}{3} \approx 0.7454$, 
so that $\iota (Y) \ge 1.366$. Hence $K = 1$ is a doubling characteristic distance.
Consequently ${\mathcal T}(F_3) \le 2 \vert B^{F_3}(1) \vert = 14$.
As $F_3$ is a subgroup of $F_2$ one may improve the previous estimate to
$$
   {\mathcal T}(F_2) \, \le \, 14
$$
by Observation 19. \smallskip

   (ii) Consider again the Cayley graph $X$ of $F_2$. Its isoperimetric constant
is  precisely  $\iota (X) = deg(X) - 2 = 2$ by Example 47. Hence
$K_X   =  \left\lceil \frac{ \log 2 }{\log 3 } \right\rceil =  1$ 
is a doubling characteristic distance; thus
$$
  {\mathcal T}(F_2) \, \le \, 2 \vert B^{F_2}(1) \vert 
    \, = \, 10,
$$
which compares better than the previous estimate with  ${\mathcal T}(F_2) = 4$. 
 \smallskip

   These computations indicate that some effort should be given to sharpen the
isoperimetric inequality of Lemma 50 used above (see Question 62.(a)).

\section*{IV.2. Tarski number for Ol'shanskii groups and for Burnside groups}

\noindent {\bf 58. On Ol'shanskii groups.}
We consider first a family of groups investigated in \cite{Ol1}.
(See also \cite{Ol2} both for this family and for other ones, discovered by the
same author, and relevant for the subject discussed here.) For any 
$\epsilon > 0$, there exists one of these  groups given as a quotient 
$\pi \colon F_2 \twoheadrightarrow G$ for which the relative growth $\alpha_G$ satisfies
$$
   \sqrt 3 \, < \, \alpha_G \, \le \, \sqrt 3 + \epsilon 
$$
and which is consequently non-amenable. Moreover Ol'shanskii has shown that
these groups do not have any non-abelian free subgroups; thus their Tarski number
satisfy ${\mathcal T}(G) \ge 5$, and ${\mathcal T}(G) \ge 6$ in case of torsion groups
(Proposition 21).  From the relative growth extimate above and from Theorem 56 (see also 
the first computation of Item 57), one obtains the
following.
\bigskip

\noindent {\bf 59. Proposition.} {\it There exists a two-generator non-amenable 
torsion-free  group $G$ without non-abelian free subgroup, 
for which the Tarski number ${\mathcal T}(G)$ satisfies 
$$
   5 \, \le \, {\mathcal T}(G) \, \le \, 34.
$$
There exist a two-generator non-amenable torsion group $G$, with all proper
subgroups cyclic, for which
$$
   6 \, \le \, {\mathcal T}(G) \, \le \, 34 .
$$
(The constructions of these groups are due to Ol'shanskii.)}
\bigskip

\noindent {\bf 60. On Burnside groups.}
   We consider next the Burnside group $B(m,n)$, given as the quotient of the
free group $F_m$ of rank $m \ge 2$ by the normal subgroup generated by
$\left\{ x^n \right\} _{x \in F_m}$, for an odd integer $n \ge 665$. It is
obvious that $B(m,n)$ does not contain any free group not reduced to $\{1\}$.
It is known that $B(m,n)$ is infinite, indeed of exponential growth (see 
VI.2.16 in \cite{Ady1}),  and indeed  not amenable \cite{Ady2}. \par

   From Theorem 3 and the last but one line in\footnote{
There are printing mistakes in the English version of \cite{Ady2}. In 
Theorem 3 of this paper, first the $C$ should read $G$, and second the
exponent of $(2m-1)$ should read 
$$
   \left[ \frac{1}{2} + \frac{ \beta }{ \gamma_R } + \frac{4}{ \delta_R } 
   \Bigg( 
           log_{2m-1} 
           \bigg( e \Big( 1 + \frac{\delta_R}{4 \gamma_R} \Big) \bigg)
    \Bigg) 
   \right]  
$$
(with the largest parenthesis $()$ as above). 
Also, in the last but one line of the paper, $\frac{1}{15}+ \frac{6}{57}$
should be replaced by $\frac{1}{15}+ \frac{5.69}{57}$, which is indeed a
number strictly smaller than $\frac{1}{6}$!
} \cite{Ady2},
one has the relative growth estimate
$$
   \alpha \, \le \, (2m-1)^{\frac{1}{2} + \frac{1}{15}+ \frac{5.69}{57}}
$$
where $\frac{1}{2} + \frac{1}{15}+ \frac{5.69}{57}$ is strictly smaller than, 
but near, $\frac{2}{3}$. \par

   For $m = 2$, Theorem 56 shows that one has successively 
$ \alpha < \root 3 \of { 9 }$, hence
$\rho   < \frac{\sqrt 3}{4} \left( 
\frac{ \sqrt 3}{\root 3 \of {9}} + \frac{ \root 3 \of {9} }{ \sqrt 3} 
\right)
\approx  0.881$, hence
$\iota (X)  \ge  4 \frac{ 1 - \rho(X) }{\rho(X)} \approx 0.540$, hence
$K = \, \left\lceil \frac{ \log 2 }{\log ( 1.540 )} \right\rceil =  2$, 
hence finally
$$
 {\mathcal T}(B(2,n)) \, \le \, 2 \vert B^{F_2} (2) \vert 
       \, = \, 2 \big( 2^.3^2 - 1 \big)
       \, = \, 34.
$$
For $m = 3$, the corresponding computations give
$ \alpha < \root 3 \of { 25 }$, hence
$\rho < \frac{\sqrt 5}{6} \left( 
\frac{ \sqrt 5}{\root 3 \of {25}} + \frac{ \root 3 \of {25} }{ \sqrt 5} 
\right)
\approx  0.772$, hence
$\iota \ge 1.181$, hence
$K = 1$, hence finally
$$
  {\mathcal T}(B(3,n)) \, \le \, 2 \vert B^{F_3} (1) \vert 
       \, = \,  14.
$$

   Let $m_1,m_2$ be such that $2 \le m_1 \le m_2 \le \infty$ and let
$n$ be as above. It follows from general principles 
on relatively free groups in varieties of groups that
$B(m_2,n)$ has a subgroup isomorphic to $B(m_1,n);$ see \cite{NeuH},
Statements 12.62 and 13.41. It is also known that $B(m_1,n)$ 
has a subgroup isomorphic to $B(m_2,n);$ see \cite{Sir},
and also \S \ 35.2 in \cite{Ol2}. 
Thus, it follows from Observation 10 that  one has
${\mathcal T} \big(B(m,n)\big) = {\mathcal T}\big(B(3,n)\big)$
for any $m \ge 2$. \par

   This and Proposition 21 show the following.
	\bigskip

   \noindent{\bf 61. Theorem.} {\it For $m \ge 2$ and for $n$ odd and at least $665$,  the
Tarski number  of the Burnside group $B(m,n)$ satisfies 
$$
   6 \, \le \, {\mathcal T}\big(B(m,n)\big) \, \le \, 14.
$$}
\bigskip

   Let us mention that it is unknown whether, for $n$ large, $B(m,n)$
has infinite  amenable quotients. (A question of Stepin, which is Problem 9.7
of \cite{Kou}.) Similarly, one could ask what are the Tarski numbers
of non-amenable quotients of these groups.
\bigskip

\noindent {\bf 62. Questions of continuity.} \par
   {\it Question (a):} given $\epsilon > 0$, does there exist $\delta > 0$ such
that, for any quotient group $G$ of a free group $F$ with spectral radius
satisfying $\rho(G) < \rho(F) + \delta$, one has necessarily an estimate
$\iota(G) > \iota(F) - \epsilon$ for the isoperimetric constants ? 
More generally, can one sharpen the inequality 
$\iota(X) \ge 4 \frac{1 - \rho(X) }{\rho(X)}$ of Lemma 50? \par

   {\it Question (b):} given $\delta > 0$, does there exist $\eta > 0$ such
that, for any quotient group $G$ of a free group $F$ 
with minimal growth rate satisfying $\omega(G) > \omega(F)-\eta$, one has necessarily an estimate
$\rho(G) < \rho(F) + \delta$? \par

   (For $\omega(G)$, see \cite{GriH}. If the free group $F$ above is of rank
$m$ and is considered together with a free basis, recall that 
$\omega(F) = 2m-1$, $\rho(F) = \frac{\sqrt{2m-1}}{m}$,  and $\iota(F) = 2m-2$. The
coefficients $\rho(G)$ and $\iota(G)$ are of course taken with
respect to the images in $G$ of free generators in $F$.)

   Assume the two questions above have affirmative answers; then: (i) for a
convenient group $G$ of Ol'shanskii, $\iota(G) \ge 2 - \epsilon$, and $K = 1$, 
and consequently  ${\mathcal T}(G) \le 10;$ (ii) for the Burnside groups $B(2,n)$
of Theorem 61 with $n$ large enough, one would have 
$\omega(G) \ge 3 - \epsilon$ (VI.2.16 in \cite{Ady1}), 
and $K = 1$, and consequently also 
${\mathcal T}\big(B(m,n)\big) \, = \, {\mathcal T}\big(B(2,n)\big) \, \le \, 10$ 
for any $m \ge 2$ and $n$ odd large enough.

\section*{{\bf V. Supramenability}}

\section*{V.1. Supramenability and subexponential growth}

	\noindent
	{\bf 63. Definition.} {\it A pseudogroup $({\mathcal G} ,X)$ is supramenable if 
the pseudogroup $(\mathcal G_{(A)}, A )$ 
defined in Example 2.(iv) is amenable for any nonempty subset $A$ of $X$.}
\bigskip

    In case of a pseudogroup ${\mathcal W}(X)$,  Remark 3.(vi) shows that one may
read this definition in two ways. More precisely, a locally finite metric
space $X$ is supramenable if, for any subspace $A$ of $X$, one has\par
    (i) the metric space $A$ is amenable, i.e. the pseudogroup ${\mathcal W}(A)$ is
amenable, \par\noindent
or equivalently \par
   (ii) the restriction ${\mathcal W}(X)_{(A)}$ of the pseudogroup ${\mathcal W}(X)$ to $A$ 
is amenable.
\par\noindent
   Observe that supramenability of uniformly locally finite metric spaces is invariant by
quasi-isometry, because of Proposition 38. \medskip

   A finitely generated group is {\it supramenable} if it so as a metric space,
for the combinatorial distance on its Cayley graph with respect to a finite
generating set (this definition of supramenability does not depend on the choice 
of the generating set). \par
   This notion, due to Rosenblatt \cite{Ros2}, carries over to not necessarily
finitely generated groups, and indeed to  topological groups, 
but we will not use this below.  
\bigskip

 \noindent  {\bf 64. Definition.} {\it Let $X$ be a locally finite metric space; for a point
$x \in X$ and a number $r \ge 0$, we denote by $\beta ^X_x(r)$
the cardinality of the closed ball of radius $r$ around $x$ in $X$.
The space $X$ is of 
$$
\begin{aligned}
&\text{ {\it subexponential growth } if  } \qquad&
   &\limsup_{r \to \infty} \root r \of { \beta ^X_x(r) } \, = \, 1 
\hskip4truecm \cr
&\text{ {\it exponential growth } if  } \qquad&
   1 \, < \,
   &\limsup_{r \to \infty} \root r \of { \beta ^X_x(r) } \, < \, \infty \cr
&\text{ {\it superexponential growth } if } \qquad&
   &\limsup_{r \to \infty} \root r \of { \beta ^X_x(r) } \, = \, \infty .\cr
\end{aligned}
$$}
\medskip

Observe\footnote{Unlike in some other places of this paper (such as Proof 36), we insist here
that the distance between two points of $X$ is always {\it finite}.} 
that any of these holds for some $x \in X$ if and only if it holds for all
$x \in X$, and also if and only if it holds for any pair $(X',x')$ with $X'$
quasi-isometric to $X$. In particular, subexponential growth and exponential 
growth make sense for finitely generated groups, without any mention of
a generating set. \medskip

\noindent
{\bf 65. Lemma.} {\it Inside a locally finite metric space of subexponential growth, any 
subspace is also of subexponential growth.} 
\medskip

\noindent{\it Proof.} For a subspace $Y$ of a space $X$, one may choose in the
previous definition the point $x$ inside $Y$. Then the lemma follows from
the obvious inequality $\beta^Y_x(r) \le \beta ^X_x(r)$, for all $r \ge 0$. \hfill \qed
\bigskip

   For historical perspective, let us recall that a simple argument going
back to \cite{AdVS} shows that a finitely generated group which is of
subexponential growth is amenable, and indeed supramenable (Theorem 4.6 in
\cite{Ros2}). \par
   As a consequence, one has $\iota (X) = 0$ for any Cayley graph $X$ of a
finitely generated group of subexponential growth. There are further
connections between growth and isoperimetry, due to Varopoulos and others.
More precisely, consider for example a finitely generated group $G$
generated by a finite set $S$, the corresponding growth function $\beta^G_S$
defined by 
$$
  \beta^G_S (n) \, = \, \vert \left\{ \, g \in G \, \vert \,
     \text{ the $S$-word length of $g$ is at most $n$ } \, \right\}\vert
$$
for all $n \ge 0$, and the {\it isoperimetric profile} $I^G_S$ defined by
$$
   I^G_S (n) \, = \, \max_{m \le n} \quad
        \min_{F \subset X^0 \, , \,  \vert F \vert = m} \,\,\,
        \vert \partial F \vert
$$
for all $n \ge 1$, where $X^0$ denotes the vertex set of the Cayley graph
of $G$ with respect to $S$ (namely $X^0 = G$!); then, for various classes
of groups, there are quite precise estimates relating the growth function
$\beta^G_S$ and the isoperimetric profile $I^G_S;$ see in particular
\cite{CoSa} and \cite{PiSa}. \par
   In our context, the argument of \cite{AdVS} provides the following result.
	
	\bigskip

\noindent{\bf 66. Theorem.} {\it A locally finite metric space of subexponential growth is
supramenable.}
\smallskip

\noindent{Proof.} Let $X$ be a locally finite metric space of subexponential growth.
By the previous lemma, it is enough to show that $X$ is amenable; we will show
that $X$ satisfies the F\o lner condition. \par

   Consider a finite subset ${\mathcal R}$ in the pseudogroup ${\mathcal W}(X)$,
a point $x_0 \in X$ and a number $\epsilon > 0$. Set
$$
   C \, = \, \left\lceil
      \max_{ \rho \in {\mathcal R}\cup {\mathcal R}^{-1}} \sup_{x \in \alpha(\rho)}
      d ( x , \rho(x) ) \right\rceil .
$$
As
$$ 
   \limsup_{r \to \infty} \root r \of { \beta^X_{x_0}(r) } \, = \, 1
$$
there exists a strictly increasing sequence of integers 
$\left( r_k \right)_{k \ge 1}$ such that
$$
   \lim_{k \to \infty} \frac{ \beta^X_{x_0}(r_k+C) }{ \beta ^X_{x_0}(r_k) }
   \, = \, 1.
$$
Set
$$
   F_k \, = \, \text{ball of radius} \quad r_k \quad \text{centered at}
       \quad x_0 \quad \text{in} \quad X
$$
for all $k \ge 1$. \par
   As $\partial_{\mathcal R} F_k \subset {\mathcal N}_C F_k \smallsetminus F_k$ 
for all $k \ge 1$, one has
$$
   \lim_{k \to \infty} \frac{ \vert \partial_{\mathcal R} F_k \vert}
          { \vert F_k \vert } \, = \, 0
$$
so that $\left( F_k \right) _{k \ge 1}$ is a ``F\o lner sequence''\,
(see Definition 23), and this ends the proof. \hfill \qed

\bigskip

   The following criterium for graphs will be used in Section V.2.
Recall that a metric space $X$ is {\it long-range connected} if there
is a constant $C > 0$ such that every two points $x$ and $y$ in $X$ can
be joined by a finite chain of points
$$
   x_0 = x \, , \, x_1 \, , \, \hdots \, , \, x_n = y
$$
such that 
$$
   d(x_{i-1},x_i) \, \le \, C
$$
for all $i \in \{1 , \hdots , n\}$ (see Item 0.2-$A_2$ in \cite{Gro3}).
\bigskip

   \noindent{\bf 67. Proposition.} {\it A connected locally finite graph is supramenable
if and only if all its long-range connected subgraphs are amenable.}
\smallskip

  \noindent{\it Proof of the non-trivial implication.} Given a graph $X$ which is
{\it not} supramenable, we have to show that there exists a long-range connected
subset $Z$ of its vertex set $X^0$ which is not amenable (as a metric space, for 
the combinatorial distance of $X$). \par
   By hypothesis, there exists a subset $Y$ of $X^0$ and a mapping
$\phi \colon Y \to Y$ such that $\sup_{y \in Y}d(\phi(y),y) \le C$ for some constant 
$C \ge 0$, and such that $\vert \phi^{-1}(y) \vert \ge 2$ for all $y \in Y$.
Set $Z = {\mathcal N}_C (Y)$, and let $\left( Z_i \right)_{i \in I}$ be an enumeration
of the connected components of $Z$. For all $i \in I$, set  $Y_i = Y \cap Z_i$.
As $\phi$ is a $C$-bounded perturbation of the identity, one has
$\phi^{-1}(Y_i) \subset Z_i$,  and it follows that $\phi^{-1}(Y_i) \subset Y_i$, 
for all $i \in I$. Hence $Y_i$ is paradoxical for each $i \in I$. \hfill \qed

\section*{V.2. Examples with trees}

 \noindent  
Let ${\mathcal S}_2$ denote the free {\it semi}group on two generators. 
From the natural word length, one defines on ${\mathcal S}_2$ a metric making it a
uniformly locally finite
metric space which is of exponential growth, and indeed paradoxical.
Thus, any finitely generated group containing a subsemigroup isomorphic to
${\mathcal S}_2$ has a paradoxical subspace (the group being viewed as a metric space),
and consequently is not supramenable. \bigskip

\noindent {\bf 68. Question.} Does there exist a finitely
generated group which is amenable, not supramenable, and without
subsemigroup isomorphic to ${\mathcal S}_2$ ? 
\medskip

   This question is due to Rosenblatt, who conjectured the answer to
be negative (see \cite{Ros2}, just after Theorem 4.6 and after Corollary 4.20);
he also observed the following alternative for a finitely generated solvable group:
either the group has a nilpotent subgroup of finite index, and then the group is
supramenable, or the group contains ${\mathcal S}_2$ as a subsemigroup, and then the group
is not supramenable (Theorems 4.7 and 4.12 in \cite{Ros2}). \par
    However, Question 68 has been answered positively by the second
author as follows. \bigskip

   \noindent{\bf 69. Examples \cite{Gri4}.} {\it For each prime $p$, there exist 
uncountably many finitely generated $p$-groups which are \par
    $\bullet$ of exponential growth, \par
    $\bullet$ without any subsemigroup isomorphic to ${\mathcal S}_2$,
              \par
    $\bullet$ amenable, \par
    $\bullet$ not supramenable.}
\medskip

  \noindent{\it On the proof.} This involves wreath products\footnote{In the English translation of \cite{Gri4}, the Russian word for
``wreath product''\, has been incorrectly translated as
``amalgamated product''!}
$G = C_p \wr H$, where $C_p$ denotes a cyclic group of order $p$ and where
$H$ is one of the $p$-groups of intermediate growth constructed in
\cite{Gri2, Gri3}. 
\par
To show that $G$ is not supramenable, the idea is to
construct a paradoxical binary subtree in an appropriate Cayley graph of $G$.  
   As a torsion group, $G$ does not contain ${\mathcal S}_2$.
   The two other claims are straightforward. \hfill \qed
\bigskip

\noindent {\bf 70. Question.} Does there exist a finitely generated group
which is supramenable and of exponential growth ? \bigskip

   This question, formulated as Item 12.9.a and Problem C.12 of \cite{Wag}, is 
still open. \par

   One way to make the question more precise is recorded as Problem 16.11 in 
the {\it Kourovka Notebook} \cite{Kou}: does there exist a finitely generated
semigroup $S$ with cancellation having subexponential growth and such that the
group of left quotients $G = S^{-1}S$ has exponential growth? (The group of
quotients would exist, because the so-called ``Ore condition''\, 
holds; see for  example Sections 1.10 and 12.4 in \cite{ClPr}.) 
The point is that such a semigroup of subexponential growth is supramenable
and that a group of quotients of a supramenable semigroup is a supramenable
group.  \par
 
   Here is however a straightforward construction. \bigskip

\noindent{\bf 71. Example.} {\it There exists a locally finite metric space which is of
superexponential growth and which is supramenable.}
\medskip

\noindent{\it Proof.} Consider a sequence $\left( d_k \right)_{k \ge 0}$
of integers $\ge 2$ and a sequence $\left( h_k \right)_{k \ge 1}$
of integers $\ge 1$. Let $X$ be a rooted tree  
in which a vertex at distance $n$ of the root is of degree 
$$
\left\{
\aligned
  d_k \qquad \qquad &\text{if} \qquad n = \sum_{j=1}^k h_j \\
   2 \qquad \qquad &\text{otherwise}
\endaligned
\right.
$$
(given the two sequences, this completely defines the
tree up to isomorphism). \par

   If $\liminf_{k \to \infty} h_k = \infty$, a long-range connected
subspace $Y$ of the vertex set of $X$ cannot satisfy the Gromov
condition (compare with Proposition 35 above, i.e. with Corollary 4.2
of \cite{DeSS}). It follows from
Proposition 67 that $X$ is supramenable. 
\smallskip

   Now the growth sequence of $X$ with respect to the root, say $x_0$, satisfies
$$
   \beta^X_{x_0}(n+1) \, \ge \, \prod_{j=0}^{k} d_j \qquad \text{for} \qquad
     n \, = \, \sum_{j=1}^{k} h_j ,
$$
so that, if  the sequence $\left( d_k \right)_{k \ge 0}$ is increasing
rapidly enough, one has 
$$
   \limsup_{m \to \infty} \root m \of { \beta^X_{x_0}(m) } \, = \, \infty
$$
and $X$ is of superexponential growth. For example, 
if $d_j = \left( \sum_{i=1}^j h_i \right)!$,
then
$$
   \beta^X_{x_0}(n+1) \, \ge \, d_k \, = \, n! 
$$
whenever $n = \sum_{j=1}^k h_j$, and this implies
$\limsup_{m \to \infty} \root m \of { \beta^X_{x_0}(m) } \, = \, \infty$
by Stirling's formula. \hfill \qed

\bigskip

\noindent {\bf 72. Variation on the previous example.}
{\it There exists a graph of bounded degree which is of exponential growth and 
which is supramenable.} 
\medskip

\noindent{\it Proof.} Consider a rooted tree $X$ in which a vertex at distance $n$
of the root is of degree \par
   $2$ if $(k-1)k \le n < k^2$ for some $k \ge 1$, \par
   $3$ if $k^2 \le n < k(k+1)$ for some $k \ge 1$. \par
\noindent The growth function of $X$ with respect to the root satisfies
$$
    2^{\frac{k(k+1)}{2}} \, \le \, 
   \beta^X \big( k(k+1) \big) \, \le \, 3^{\frac{k(k+1)}{2}}
$$
for all $k \ge 1$, so that $X$ is clearly of exponential growth. Example 48
implies that $X$ is supramenable. \hfill \qed
\bigskip

\noindent {\bf 73. Question.} Let $G$ and $H$ be two finitely generated groups
which are supramenable; is the product $G \times H$ supramenable ?
\bigskip

   This question appears in \cite{Ros2} (just before Proposition 4.21), and the
answer is still unknown. Here is however an example, for which we are
grateful to Laurent Bartholdi. \bigskip

\noindent {\bf 74. Example.}  {\it There exist two supramenable locally finite metric spaces
$X,Y$ such that the direct product $X \times Y$ is not supramenable, for
the metric defined by
$$
   d_{X \times Y}\big( (x_1,y_1) , (x_2,y_2) \big) \, = \, 
        d_{X}(x_1,x_2) + d_{Y}(y_1,y_2).
$$ }
\medskip

\noindent {\it Proof.} Let $\left( h_k \right)_{k \ge 1}$ be a strictly
increasing sequence of integers $\ge 1$.  Let $X$ be a rooted tree 
in which a vertex at distance $n$ of the root
is of degree 
$$
\left\{
\aligned
  3 \qquad  &\text{if} \qquad  
       \sum_{j=0}^{2k} h_j \, \le \, n \, < \, \sum_{j=1}^{2k+1} h_j
        \qquad \text{for some} \quad k \ge 0, \\
   2 \qquad  &\text{otherwise}
\endaligned
\right.
$$
(with $\sum_{j=0}^{2k} h_j = 0$ for $k = 0$).
And let $Y$ be a rooted tree  in which a vertex at distance $n$ of the root
is of degree 
$$
\left\{
\aligned
  3 \qquad  &\text{if} \qquad  
       \sum_{j=1}^{2k+1} h_j \, \le \, n \, < \, \sum_{j=1}^{2k+2} h_j
        \qquad \text{for some} \quad k \ge 0, \\
   2 \qquad  &\text{otherwise.}
\endaligned
\right.
$$
Observe that both $X$ and $Y$ are supramenable, because each of their
infinite connected subgraphs has
arbitrarily large hanging chains. Observe also that, for each integer $n$,
there is either in $X$ or in $Y$ a vertex of degree 3 at distance $n$ of the
relevant root. It follows that the product of the two
metric spaces defined by $X$ and $Y$, for the distance 
$d_{X \times Y}$ defined above,
contains a paradoxical tree. Consequently, $X \times Y$ is not
supramenable. \hfill \qed

\bigskip

\noindent {\bf 75. Paradoxical subtrees in paradoxical graphs.}
   It is known that a paradoxical graph contains a paradoxical tree
\cite{BeSc}. It is unknown whether a connected paradoxical graph necessarily
contains a paradoxical tree which is {\it spanning,} i.e. which contains all 
vertices of the original graph (this is Problem 2 in Section 4 of \cite{DeSS}). \par
   However, Benjamini and Schramm have shown that, if $X$ is a paradoxical
graph with  $\iota(X) \ge n$ for some integer $n \ge 2$, 
then $X$ has a {\it spanning forest} of which every connected component is a
tree with one vertex of degree $n - 1$ and all other vertices of degree $n + 1$.  
This implies that $X$ has a paradoxical {\it spanning tree.} 

\bigskip

\noindent {\bf 76. A question of V. Trofimov.} This appears as Problem
12.87 in the Kourovka Notebook \cite{Kou}.
   Let $X$ be a connected  undirected graph without loops and multiple edges
and suppose that its automorphism group $Aut(X)$ acts transitively on the
vertices. Is it true that one of the following holds? \par

   (i) the stabilizer of a vertex of $X$ is finite, \par
   (ii) the action of $Aut(X)$ on the vertices of $X$ admits a
non-trivial imprimitivity system $\sigma$ with finite blocks for which the
stabilizer of a vertex of the factor-graph $X/\sigma$ in $Aut(X /\sigma)$
is finite, \par
   (iii) there exists a natural number $n$ such that the graph,
obtained from $X$ by adding edges connecting distinct vertices the distance
between which in $X$ is at most $n$, contains a tree all of whose vertices
have valence $3$.
\par

  If the answer to this question was positive, this would imply that a graph of
subexponential growth having a transitive group of automorphisms is essentially
a Cayley graph of a group.

\section*{{\bf VI. Comments and corrections (March 2016)}}

\noindent \textbf{On the article by Deuber, Simonovits and S\' os, 
and their terminology of exponential growth.}
There is an annotated version of the 1995 version, dated 2004 \cite{DeSS--04},
and an exposition of related material \cite{ElSo--05}. 
In the annotated version, the authors observe that their terminology
of exponential growth is not standard in the group theory literature.

\vskip.2cm
\noindent \textbf{On No.\ 11 and elementary amenable groups.}
Let $B_0$ denote the class consisting of all finite groups and the infinite cyclic group.
The following fact was shown by Chou and refined by Osin \cite[Theorem 2.1]{Osin--02}:
the class $EG$ of elementary amenable groups is the smallest class of groups
which contains the trivial group $\{1\}$, 
which is closed under taking direct limits,
and which is such that a group $G$ is in $EG$ whenever there exists an extension
$\{1\} \to N \to G \to Q \to \{1\}$ with $N \in EG$ and $Q \in B_0$.
\par
Nekrashevych \cite{Nekr} has recently discovered examples 
of finitely generated infinite groups
that are simple, periodic, and of intermediate growth. 
In particular
they are amenable (because of intermediate growth)
and not elementary amenable (either because 
infinite finitely generated simple groups cannot be elementary amenable,
or because infinite finitely generated periodic groups cannot be elementary amenable --
both observations go back to Chou).

\vskip.2cm
\noindent \textbf{On No.\ 13 and the space of marked groups.}
The space of marked groups has received
a considerable amount of attention.  
Besides the articles cited in No.\ 13, we indicate
\cite{ChGu--05}, \cite{CoGP--07}, and \cite{BCGS--14}.
\par
See also \cite[Corollary 6.25]{WeWi} for an original proof of the existence of finitely generated groups
that are amenable and are not elementary amenable:
in the appropriate space of marked groups, the set of amenable groups is Borel
and the set of elementary amenable groups is not.

\vskip.2cm
\noindent \textbf{On No.\ 14 and subexponentially amenable groups.}
There are amenable groups (i.e.\ groups in $AG$) 
that are not subexponentially amenable (i.e.\ not in $BG$). 
Indeed, the so-called \emph{Basilica group}
was first shown to be not in $BG$ \cite{GrZu--02},
and later shown to be amenable \cite{BaVi--05}. 
The method of Bartholdi and Virag 
was streamlined and generalized in \cite{Kaim--05}.
Further examples can be found in
\cite{Ersc--06}, \cite{Brie--09}, \cite{BaKN--10}.
The finitely generated amenable simple groups
that appear in \cite{JuMo--13} are also amenable and not subexponentially amenable.

\vskip.2cm
\noindent \textbf{On Nos.\ 15 and 24, Ahlfors' notion of regular exhaustion,
and Bogolyubov's ideas on amenability for topological groups.}
In \cite{Roe--88},
there is a discussion of regular exhaustion, introduced by Ahlfors in 1935 [Ahl].
In \cite{GrHa},
there is a discussion of Bogolyubov's ideas on  amenability, in his 1939 article [BogL]
which went almost unnoticed.
\par
There is an exposition of basic material on amenability 
of topological groups (and the important case of locally compact groups)
in Chapter II.G of \cite{BeHV--08}.

\vskip.2cm
\noindent \textbf{On No.\ 15, amenability of groups and cellular automata.}
Let $G$ be a group and $A$ a finite set.
Equip $A^G = \{u \colon G \to A\}$ with its \emph{prodiscrete topology} (i.e.~the topology of pointwise convergence) 
and with the \emph{shift action} of $G$ defined by $gu(h) := u(g^{-1}h)$ for all $g,h \in G$ and $u \in A^G$.
A \emph{cellular automaton} over $G$ is a continuous map $\tau \colon A^G \to A^G$ that is $G$-equivariant, i.e., satisfies $\tau(g u) = g \tau(u)$ 
for all $g \in G$ and $u \in A^G$.  
For two maps $u, v  \in A^G$ write $u \approx v$ if they coincide outside of a finite subset of $G$. It is clear that $\approx$ is an equivalence relation.
A map $\tau \colon A^G \to A^G$ is said to be  \emph{pre-injective} if its restriction to each $\approx$-equivalence class is injective.
Then the \emph{Garden of Eden theorem}~\cite{CeSMS--99} (see also \cite{Grom--99b}), originally established by Moore \cite{Moo--63} and Myhill \cite{Myh--63} 
for $G = \Z$,  states that a cellular automaton over an amenable group is surjective if and only if it is pre-injective.
It follows from\cite{Bart--10} that if a group $G$ is non-amenable then there exist cellular automata over $G$ that are surjective but not pre-injective. Thus, the Garden of Eden theorem yields a characterization of amenability for groups in terms of cellular automata. 
For more on the Garden of Eden theorem (consequences and variations) we refer to \cite{CeSC--10}.


\vskip.2cm
\noindent \textbf{On No.\ 17, inner amenability, and coamenability.}
An old question on inner amenability from \cite{Eff} has been solved in \cite{Vaes--12}.
\par
Several claims of Theorem 5 in [BeHa] have to be corrected,
as in \cite[Section 3]{Stal--06}.
\par
For the notion of coamenability of a subgroup of a group,
see also \cite{MonPo--03} and \cite{Pest--03}.

\vskip.2cm
\noindent \textbf{On Definition 18, Question 22,  and Tarski numbers.}
There are other definitions of Tarski numbers, see \cite[Appendix A]{ErSP--15}.
Since 1999, there has been some progress on understanding of Tarski numbers.
For example, 
there are $2$-generated non-amenable groups  with arbitrarily large Tarski numbers, 
there are groups which we know have Tarski number exactly $5$, or $6$,
and every number $\tau \ge 4$ is the Tarski number of some
faithful transitive action of a finitely generated free group.
See \cite{OzSa--13}, \cite{ErSP--15}, \cite{Gola--a}, and \cite{Gola--b}.

\vskip.2cm
\noindent \textbf{On Definition 29 and the reference [GrLP].}
In its second edition, this book has been considerably expanded \cite{Grom--99a}.

\vskip.2cm
\noindent \textbf{On Definition 30 and the terminology ``doubling condition''.}
It is unfortunate that this terminology is used in several incompatible meanings.
Some authors use them in our sense, see e.g.\ \cite{Kapo--02}.
But many more authors use them in a completely different meaning,
most often for metric spaces with measures, 
and occasionally for metric spaces as such;
see e.g.\ \cite{Grom--99a}, \cite{Hein--01} and \cite{LoVi--07}.
\par
More precisely, a metric space $X$ is called \emph{doubling} if there exists a constant $C > 0$
such that, for all $d > 0$, any subset of $X$ of diameter at most $d$
can be covered by $C$ subsets of $X$ of diametrer at most $d/2$ \cite[Definition 10.13]{Hein--01}.
The doubling metric spaces are precisely the spaces of
\emph{finite Assouad dimension}; compare \cite[Definition 10.15]{Hein--01}.
\par

In retrospect, our terminology for the notion of Defintion 30 was unfortunate.
\par

A change of terminology there should have some effect on the terminology
``doubling characteristic distance'' of No.\ 53.

\vskip.2cm
\noindent \textbf{On Section III.1, discrete and locally finite metric spaces,
and uniform notions.}
In the published version, just before Defintion 28,
we have unfortunately used the word ``discrete''
for what should be ``locally finite''.

For a metric space $(X,d)$, the four following properties should not be confused:
\begin{enumerate}[noitemsep]
\item[$\bullet$]
$(X,d)$ is \emph{discrete} if, for every $x \in X$, there exists $\delta_x > 0$
such that $d(x,y) \ge \delta_x$ for all $y \in X \smallsetminus \{x\}$;
note that $(X,d)$ is discrete if and only if the topology on $X$ defined by $d$ is discrete;
\item[$\bullet$]
$(X,d)$ is \emph{uniformly discrete} if there exists $\delta > 0$
such that $d(x,y) \ge \delta$ for all $x,y \in X$ such that $x \ne y$;
\item[$\bullet$]
$(X,d)$ is \emph{locally finite} if  every subset of $X$ of finite diameter is finite;
\item[$\bullet$]
$(X,d)$ is \emph{uniformly locally finite} if, for every $D \ge 0$,
there exists a constant $C$ such that every subset of $X$ of diameter at most $D$
has at most $C$ elements.
\end{enumerate}
Note that a discrete metric space is locally finite if and only if it is proper,
i.e.\ if and only if its closed balls are compact.
Note also that a uniformly locally finite metric space need not be uniformly discrete
(example: the subspace $\{ n \in \Z \mid n \ge 1\} \cup \{ n+2^{-n} \mid n \ge 1 \}$
of the real line).
Let $X$ be a connected graph, and $(X^0, d)$ its vertex set
together with the combinatorial distance function;
then $X^0$ is always uniformly discrete,
and $X^0$ is locally finite [respectively uniformly locally finite]
if and only if $X$ is locally finite [respectively of bounded degree].
\par

In the present version (unlike in the published version), 
we have used ``locally finite'' instead of ``discrete''
in Nos.\ 28 to 36.
Proposition 38, on invariance of amenability by quasi-isometries,
holds for \emph{uniformly} locally finite metric spaces.

\vskip.2cm
\noindent\textbf{An example in \cite{DiMW}.}
Here is the last example of \cite{DiMW}, showing that
the hypothesis of uniform local finiteness cannot be deleted in Proposition 38.
\par

Consider the graph $X$ defined as follows: 
\begin{enumerate}[noitemsep]
\item[$\cdot$]
it has vertices $(n,1)$ for all $n \in \Z$ with $n \ge -1$,
\item[]
and $(n,k)$ for all $n \ge 1$ and $k \in \N$ such that $2 \le k \le 2^n$;
\item[$\cdot$]
it has  edges connecting $(n,1)$ to $(n+1,1)$ for all $n \in \Z$ with  $n \ge -1$,
\item[]
and $(n,1)$ to $(n,k)$ for all $n \ge 1$ 
and $k \in \N$ such that $2 \le k \le 2^n$.
\end{enumerate}
Observe that $X$ is locally finite and not uniformly locally finite.
Denote by $X^0$ the vertex set of this graph,
considered as a metric space for the combinatorial metric, say $d$; 
observe that $X^0$ is a uniformly discrete metric space.
Let $\Phi : X^0 \longmapsto X^0$ be the mapping defined as follows:
\begin{enumerate}[noitemsep]
\item[$\cdot$]
$\Phi(-1,1) = \Phi(0,1) = (-1,1)$;
\item[$\cdot$]
$\Phi(n,k) = \Phi(n, k+2^{n-1}) = (n-1,k)$
for all $n \ge 1$ and $k$ with $1 \le k \le 2^{n-1}$.
\end{enumerate}
Then $\Phi$ is a bounded perturbation of the identity
and all its fibers have two elements, in other terms
$$
d(\Phi(x), x) \, \le \, 3
\hskip.2cm \text{and} \hskip.2cm
\vert \Phi^{-1}(x) \vert = 2
\hskip.2cm \text{for all} \hskip.2cm x \in X^0 .
$$
Hence $X^0$ satisfies the Gromov condition of Definition 29,
and $X^0$ is paradoxical by Theorem 32.
Consider also the subgraph of this graph with vertex set
$Y^0 = \{ (n,1) \mid n \in \Z, n \ge -1 \}$
and edges connecting $(n,1)$ to $(n+1,1)$ for all $n \in \Z$ with $n \ge -1$;
observe that this graph is a half line, 
and that the corresponding discrete metric space $Y$ is amenable.
\par
There is an obvious quasi-isometry from $X^0$ to $Y^0$,
that maps $(n,k)$ to $(n,1)$ for all $(n,k) \in X^0$.
Yet $X^0$ is paradoxical and $Y^0$ amenable.

\vskip.2cm
\noindent \textbf{On metric spaces
for which amenability could make sense.}
Logically, Definitions 28, 29, 30 would make sense for every metric space.
In Theorem 32 implications 
$$
\begin{array}{*{20}c}
 & & & & (v) &&&& (vi) 
 \\
 &&&&&&&&
 \\
&&&& \Updownarrow &&&& \Downarrow 
 \\
&&&&&&&&
 \\
 (vi) & \Leftarrow & (i) & \Leftrightarrow & (ii) & \Rightarrow & (iii) & \Rightarrow & (iv)
\end{array}
$$
would still be correct.
But local finiteness is important for our proof of (iv) $\Rightarrow$ (v).
Indeed, Hall-Rado Theorem, No.\  35, does not carry over to arbitrary bipartite graphs,
as the following example shows.
\par
Consider the bipartite graph $B = B(Y,Z; E)$ 
of which the vertex set is the disjoint union of two sets
given with bijections with the integers, 
say $\alpha : Y \longrightarrow \N$ and $\beta : Z \longrightarrow \N$ 
(recall that $\N$ contains $0$),
and the edge set is
$$
E \, = \, \{ (y,z) \in Y \times Z \mid \alpha(y) = \beta(z)+1 \}
\sqcup
\{ (y,z) \in Y \times Z \mid \alpha(y) = 0 \} .
$$
in other terms, $\alpha^{-1}(n+1) \in Y$ has a unique neighbour 
$\beta^{-1}(n)$ for all $n \ge 0$,
and the set of neighbours of $\alpha^{-1}(0)$ is the whole of $Z$.
Then $\vert \partial_E F \vert \ge \vert F \vert$ for all every finite subset $F$
of either $Y$ or $Z$,
but there does not exist any $(1,1)$-matching of~$B$,
i.e.\ $B$ satisfies the hypothesis of Hall-Rado Theorem, but not the conclusion.
\par

We wish to stress that local finiteness is important for our Theorem 32,
and than an even stronger condition, 
\emph{uniform} local finiteness, is important for Proposition 38.

\vskip.2cm
\noindent \textbf{On Remark 42 and metric spaces
for which amenability does make sense.}
A subspace $Y$ of a metric space $X$ is \emph{cobounded} if
$\sup_{x \in X} d(x,Y) < \infty$. 
A subspace $Y$ of $X$ which is both uniformly discrete and cobounded
is a \emph{net}, as defined in Remark 42, also called 
a \emph{metric lattice} \cite[Section 3.C]{CoHa}.
An application of Zorn Lemma shows that
every metric space contains metric lattices.
A metric space $X$ is \emph{uniformly coarsely proper}
if there exists $R_0 \ge 0$ such that, for every $R \ge 0$,
there exists an integer $N$ such that
every ball of radius $R$ in $X$ can be covered by $N$ balls of radius $R_0$,
equivalently if $X$ contains a uniformly locally finite metric lattice
(for this equivalence, and others, see \cite[Proposition 3.D.16]{CoHa}).
\par
For uniformly coarsely metric spaces, 
amenability makes good sense, 
and is invariant by coarse equivalence,
in particular is invariant by quasi-isometry.

\vskip.2cm
\noindent \textbf{On Lemma 50 and an isoperimetric inequality.}
A better inequality than that of the end of No.\ 50 appears in
\cite[Theorem 3.1(b)]{Moha--88}.
Particularized to our situation (regular graph) and with our notation,
it reads
$$
\iota_*(X) \, \ge \,  \frac{d^2}{d-1-(1-\rho(X))} (1 - \rho(X)) 
$$
(note that $\rho(X) \le 1$). 
We are grateful to T.\ Nagnibeda for this reference.

\vskip.2cm
\noindent
\textbf{On No.\ 52 and the formula expressing $\rho$ in terms of $\alpha$.}
For an elaboration of this formula, see \cite{Bart--99}.
There, Bartholdi establishes an equality between two generating functions,
one related to numbers of circuits of length $n$ in some appropriate graph,
and the other related to numbers of circuits of length $n$ with no backtracking in the same graph;
the formula of No.\ 52 is then obtained as the equality between the
radii of convergence of these two formal power series.

\vskip.2cm
\noindent \textbf{On No.\ 61 and infinite amenable quotients of Burnside groups.}
The existence of such quotients still appears as an open question
in a more recent edition of the Kourovka Notebook \cite{Kour--15}.

\vskip.2cm
\noindent \textbf{On Question 62.b and amenable quotients of $F_m$
having large growth rate.}
Let $G$ be a finitely generated group and $S$ a finite generating set.
For every integer $k \ge 0$, denote by $\beta^G_S(k)$ the number of
elements $g \in G$ that can be written 
as products of at most $k$ elements in $S \cup S^{-1}$.
The \emph{exponential growth rate} of the pair $(G,S)$ is the limit
$\omega(G,S) = \lim_{k \to \infty} \root k \of { \beta^G_S(k) }$;
the existence of the limit follows from the submultiplicativity of the sequence
$(\beta^G_S(k))_{k \ge 0}$.
\par

The answer to the analogue for exponential growth rates 
of Question (b) in No.\ 62 is negative; indeed
the following is shown in \cite{ArGG--05}.
Consider an integer $m \ge 2$ and the free group $F_m$ of rank $m$;
there exists a sequence $(N_n)_{n \ge 1}$ of normal subgroups of $F_m$
such that the quotient group $G_n := F_m/N_n$ is amenable for all $n \ge 1$
and $\lim_{n \to \infty} \omega(G_n, S_n) = 2m-1$,
where $S_n$ stands for the image of a free generating set of $F_m$
by the canonical projection of $F_m$ onto $G_n$.
Moreover, the sequence $(N_n)_{n \ge 1}$ can be chosen such that
$G_n$ is abelian-by-nilpotent for all $n \ge 1$, 
or metabelian-by-finite for all $n \ge 1$.
\par

\vskip.2cm

The \emph{minimal growth rate} of a finitely generated group $G$
is the number $\omega(G) = \inf_S \omega(G,S)$,
where the infimum is taken over all finite generating sets $S$ of $G$.
For a group which can be generated by $m$ elements,
it is standard that $1 \le \omega(G) \le 2m-1$,
with equality on the right if and only if $G$ is free of rank $m$.
For this and more on minimal growth rates, see \cite{GriH}.
\par
As much as we know, Question 62.b itself, on $\omega (G)$, is still open:
For $m \ge 2$, does there exist a sequence $(N_n)_{n \ge 1}$
of normal subgroups of $F_m$ such that the quotient group $G_n = F_m/N_n$
is amenable for all $n \ge 1$ and $\lim_{n \to \infty} \omega(G_n) = 2m-1$?

\vskip.2cm
\noindent
\textbf{A misprint in No.\ 62.}
Watch out:
with the normalization chosen in our paper 
(the same as in the original paper \cite{Kest--59}),
$\rho(F_m) = \sqrt{2m-1}/m$; the value $\sqrt{2m-1}$
of the published version of our paper is a misprint.
\par

\vskip.2cm

\noindent
\textbf{On No.\ 63, and equivalent definitions of supramenability for groups.}
The following is established (among other things) in \cite{KeMR--13}.
For a group $G$ (which need not be finitely generated), the following conditions are equivalent:
\begin{enumerate}
\item[$\bullet$]
$G$ is supramenable,
\item[$\bullet$]
every cocompact action of $G$ on a locally compact Hausdorff space
admits a non-zero invariant Radon measure,
\item[$\bullet$]
there is no injective Lipschitz map from the free group of rank two to $G$.
\end{enumerate}
A map $f$ from a group $G$ to a group $H$ is \emph{Lipschitz} if,
for every finite subset $S$ of $G$, there exists a finite subset $T$ of $H$ such that
$f(x)f(y)^{-1} \in T$ for every $x,y \in G$ with $xy^{-1} \in S$.

\vskip.2cm

\noindent
\textbf{On No.\ 64, and types of growth for locally finite metric spaces.}
The notions of this definition, i.e.\ subexponential growth, exponential growth,
and superexponential growth, should be restricted to
uniformly locally finite metric spaces 
(rather than to discrete metric spaces as in the 1999 publication).
\par
Note that they are meaningfull for locally finite metric spaces,
but in this context they are not invariant by quasi-isometry.
See the example given above in the comment on Proposition 38,
or \cite[Example 3.D.7]{CoHa}.

\vskip.2cm

\noindent
\textbf{On isoperimetric profiles, as in the remark that follows Lemma 65.}
For precise estimates of various isoperimetric profiles -- or, equivalently, of the corresponding
\emph{F\o lner functions} -- see \cite{Ersc--03}.

\vskip.2cm

\noindent
\textbf{On the terminology used for Proposition 67: coarsely connected metric spaces.}
Rather than ``long-range connected'',
here is the terminology used in various places, including \cite{CoHa}:
a metric space $X$ is \emph{coarsely connected} if there exists a constant $C > 0$ such that
for every pair of points $(x,x')$ in $X$, there exists a finite sequence of points
$(x_0 = x, x_1, \hdots, x_n = x')$ in $X$ such that $d(x_{i-1}, x_i) \le C$ for all 
$i \in \{1, \hdots, n\}$. 
The point is that coarse connectedness is invariant by coarse equivalence.

\vskip.2cm

\noindent
\textbf{On Questions 70 and 73 on supramenable groups.}
At the best of our knowledge, these two questions are still open:
\begin{enumerate}[noitemsep]
\item[$\bullet$]
does there exist a supramenable group of exponential growth?
(Rosenblatt's question);
\item[$\bullet$]
 is it true that the direct product of two supramenable groups
is always supramenable?
\end{enumerate}

\vskip.2cm
\noindent \textbf{On Question 76, of Trofimov.}
This appears still as an open question
in a more recent edition of the Kourovka Notebook \cite{Kour--15}.

\end{document}